\documentclass[12pt]{article}
\usepackage{a4wide}
\usepackage{hyperref}
\usepackage{amsmath,subfigure}
\usepackage{amssymb}
\usepackage{epsfig}
\usepackage{float}
\usepackage{latexsym,graphics}
\usepackage{dsfont}
\usepackage{color}

\newtheorem{theorem}{Theorem}[section]
\newtheorem{proposition}[theorem]{Proposition}
\newtheorem{lemma}[theorem]{Lemma}
\newtheorem{corollary}[theorem]{Corollary}

\newcommand {\f}   {\frac}

\newcommand {\p}   {\partial}

\newcommand{\dis}{\displaystyle}

\def \ep{\varepsilon}
\def \R{\mathbb{R}}
\def \U{{\cal{U}}}
\def \P{{\cal{P}}}
\def \D{{\cal{D}}}
\def \G{{\cal{G}}}
\def\sqr#1#2{{\vcenter{\vbox{\hrule height.#2pt\hbox{\vrule width.#2pt height#1pt \kern#1pt\vrule width.#2pt}\hrule height.#2pt}}}}
\def\qed{\hfill$\sqr45$\bigskip}

\newcommand{\jpz}[2]{{#2}} % this command highlights the changes in blue and leaves
\title{A Structured Population Model of Cell Differentiation}

\author{Marie Doumic \thanks{{INRIA Paris-Rocquencourt, EPI BANG, 
           Domaine de Voluceau, F~78153 Le Chesnay cedex;
email: marie.doumic@inria.fr}} \footnotemark[3]
\and 
Anna Marciniak-Czochra\thanks{University of Heidelberg
Center for Modelling and Simulation in the Biosciences (BIOMS)
Interdisciplinary Center for Scientific Computing (IWR)
Institute of Applied Mathematics and
BIOQUANT}
\and Beno\^ \i t Perthame\thanks{{ Universit\'e Pierre et Marie Curie, CNRS UMR 7598 LJLL, BC187, 4, place Jussieu,  F-75252 Paris cedex 5; email: benoit.perthame@upmc.fr}}  \footnotemark[1] 
\and Jorge P. Zubelli\thanks{IMPA, Est. D. Castorina 110, Rio de Janeiro, RJ 22460-320 Brazil;
email: zubelli@impa.br}}

\date{\today}

\begin{document}
\date{\today}
\maketitle

\begin{abstract}\jpz{}{
We introduce and analyze several aspects of a new model for cell differentiation. It assumes that differentiation of progenitor cells is a continuous process. From the mathematical point of view, it is based on partial differential equations of transport type. Specifically, it consists of a structured population equation with a nonlinear feedback loop. This models the signaling process due to cytokines, which regulate the differentiation and proliferation process. 
We compare the continuous model  to its discrete counterpart,  a multi-compartmental model of a discrete collection of cell subpopulations recently proposed by Marciniak-Czochra {\it et al.}~\cite{Anna1} to investigate the  dynamics of the hematopoietic system.
We obtain uniform bounds for the solutions, characterize steady state solutions, and analyze their linearized stability. We show how persistence or extinction might occur according to values of parameters that characterize the stem cells self-renewal. 
We also perform numerical simulations and discuss the qualitative behavior of the continuous model {\it vis a vis} the discrete one. 
}

\end{abstract}


\tableofcontents

%------------------------------
\section*{Introduction}
%-------------------------------

Cell differentiation is a process by which dividing cells become specialized and equipped to perform specific functions such as nerve cell communication or muscle contraction. Differentiation occurs many times during the development of a multicellular organism as the organism changes from a single zygote to a complex system \jpz{of}{with} cells of different types. Differentiation is also a common process in adult tissues. During tissue repair and during normal cell turnover a steady supply of somatic cells is ensured by proliferation of corresponding adult stem cells, which retain the capability for self-renewal. Also various cancers are likely to originate from a population of cancer stem cells that have properties comparable to those of stem cells \cite{Al-Hajj}.  

Stem cell state and fate depends on the environment, which ensures that the critical stem cell character and activity in homeostasis is conserved, and that repair and development are accomplished \cite{Moore}.  Cell differentiation and the maintenance of self-renewal are intrinsically complex processes requiring the coordinated dynamic expression of hundreds of genes and proteins in response to external signaling. During differentiation, certain genes become activated and other genes inactivated in an intricately regulated fashion.  As a result, differentiated cells develop specific structures and performs specific functions. There exists evidence that disorder in self-renewal behavior may lead to neoplasia \cite{Al-Hajj,Beachy}. For example, it has been shown that acute myeloid leukemia originates from a hierarchy of cells that differ with respect to self-renewal capacities \cite{Hope,Schroeder}.  Although much progress has been made in identifying the specific factors and genes responsible for stem cells decisions \cite{Morrison}, the mechanisms involved in these processes remain largely unknown.
 
While different genetic and epigenetic processes are involved in formation and maintenance of different tissues, the dynamics of population depends on the relative importance of symmetric and asymmetric cell divisions, cell differentiation and death. The same genes and proteins are observed to be essential for regulation of different tissues \cite{Weissmann1}. This unity and conservation of basic processes implies that their mathematical models can apply across the spectrum of normal and pathological (cancer stem cells) development.

One established method of modeling such systems is to use a discrete collection of ordinary differential equations describing dynamics of cells at different maturation stages and transition between the stages. 
 These so called multi-compartmental models are based on the assumption that in each lineage of cell precursors there exists a discrete chain of maturation stages, which are sequentially traversed, e.g., \cite{Lord,Uchida}. However, it is also becoming progressively clear that the differentiated precursors form such sequence only under homeostatic (steady-state) conditions. Committed cells generally form a continuous sequence, which may involve incremental  stages, part of which may be reversible. As an example, cell differentiation without cell divisions is observed during neurogenesis. Moreover, in some tissues such as the mammary gland, different stages of differentiation are not well identified \cite{Dontu}. 
 
These observations invoke not only the fundamental biological question of whether cell differentiation is a discrete or a continuous process and what \jpz{the is the}{is the} measure of cell differentiation, but also how to choose an appropriate modeling approach.  Is the pace of maturation (commitment) dictated by successive divisions, or is maturation a continuous process decoupled from proliferation? In leukemias, it seems to be decoupled. The classical view in normal hematopoiesis seems to be opposite. 
 
To address these questions and to investigate the impact of possible continuous transformations on the differentiation process, we introduce a new model based on partial differential equations of transport type and compare this model to its discrete counterpart. The point of departure is a multi-compartmental model of a discrete collection of cell subpopulations, which was recently proposed in \cite{Anna1} to investigate dynamics of the hematopoietic system with cell proliferation and differentiation regulated by a nonlinear feedback loop. Furthermore, since self-renewal is an important parameter in our models, the proposed models seem to be a right departure point to investigate cancer development, for example in leukemias \cite{Schroeder}.

In the present paper we extend the discrete model to a structured population model accounting for a continuous process of differentiation of progenitor cells. Models of the latter type have been already applied to \jpz{}{the} description of some aspects of hematopoiesis  \cite{Mackey, 2Mackey, Adimy1, Adimy2, Kim}. These models are based on the assumption that differentiation of progenitor cells is a continuous process, which progresses with a constant velocity. Mathematical description involves so called age-structured population equations.
The model presented here is novel due to the nonlinearities in the coupling of the model equations, in particular the nonlinear coupling in the maturity rate function.

The paper is organized as follows. In Section \ref{Section_model}, we formulate the new model. In Section \ref{DiscCont}  the link between our model and the discrete model of \cite{Anna1} is accomplished. Sections \ref{Analysis}, \ref{Section_extinct} and \ref{Section_stationary}  are devoted to the analysis of the model. In Section \ref{Analysis}  the existence and uniform boundedness of the solutions are shown. In Section \ref{Section_extinct}, it is shown that depending on the value of a parameter characterizing stem cells self-renewal model solutions tend to zero or they stay separated from zero. Section \ref{Section_stationary} provides the structure of steady states and  conditions for existence of a positive stationary solution, while Section~\ref{Stability} is devoted to a linearized problem around the positive steady state when this state exists to investigate its stability. Using the characteristic equation we study some special cases, for which we show stability or instability of the positive stationary solution. In Section \ref{Numeric}, \jpz{}{a} numerical approach and some results on stability and instability are presented. \jpz{}{We conclude in Section~\ref{concs} with some final comments and suggestions for further investigation. }

\section{Model of Cell Differentiation}\label{Section_model}

\subsection{Continuous Model}
In the following we assume that the dynamics of differentiated precursors can be approximated by a continuous maturation model. Under this assumption  we extend the multi-compartmental system from \cite{Anna1}.  Let  $w(t)$ denote the number of  stem cells, $v(t)$  the number of mature cells and  $u(x,\, t)$ the distribution density of progenitor cells structured with respect to the maturity level $x$, so that $\int_{x_1}^{x_2}u(x,t)dx$ is equal to the number of progenitors with maturity between $x_1$ and $x_2$. This includes maturity stages between  stem cells and differentiated cells.  Thus, $u(0,t)$ describes a population of  stem cells and $u(x,\, t)$, for $x>0$,  corresponds to progenitor cells. We assume that $x= x^{*}$ denotes the last maturity level of immature cells, and therefore, $u(x^{*},t)$ describes the concentration of cells which differentiate into mature cells.

The model takes the form
\begin{eqnarray}
 \frac{d}{dt} w(t) & = & [2a_{w}(s)-1]p_{w}(s)w(t)-d_w w(t), \label{StructModel1:w}\\
\partial_t u(x,t)+\partial_x[g(x,s)u(x,t)] & = &p(x,s)u(x,t) - d(x)u(x,t),\label{StructModel1:u}\\
g(0,s)u(0,t) & = & 2[1-a_{w}(s)] p_w(s) w(t),\:\:\:t>0,\label{StructModel1:bound0}\\
\frac{d}{dt}v(t) & = & g(x^{*},s) u(x^{*},t)-\mu v(t), \label{StructModel1:v}                                              
\end{eqnarray}
together with initial data
$$
w(0)=w_0\geq0,\qquad u(0,x)=u_{0}(x) \geq 0,\qquad  v(0)=v_{0}\geq 0.
$$
Integrating formally Equation \eqref{StructModel1:u} and adding it to Equations \eqref{StructModel1:w} and \eqref{StructModel1:v} yields the following cell number balance equation
\begin{equation}\label{massbalance1}
\frac{d}{dt}\big[w(t) +\int_0^{x^*} u(x,t)dx + v(t) \big]=(p_w(s) -d_w) w(t)  +\int_0^{x^*} \big(p(x,s)-d(x)\big) u(x,t) dx - \mu v(t). 
\end{equation}
\jpz{It describes}{System~(\ref{StructModel1:w})-(\ref{StructModel1:v}) describes }the following scenario: After division a stem cell gives rise to two progeny cells. Cell divisions can be symmetric or asymmetric. We assume that on the average the fraction $a_w$ of progeny cells remains at the same stage of differentiation as the
parent cell, while the $1-a_w$ fraction of the progeny cells differentiates, i.e. transfers to the higher differentiation stage. This covers the symmetric and asymmetric scenarios. Parameters $p_w$ and $d_w$  denote the proliferation rate of stem cells and their death rate, respectively. Progenitor cells differentiate at the rate $g$, which depends on their maturity stage and is also regulated by the feedback from mature cells given by a signaling factor $s$. Parameters $p(x)$ and $d(x)$ denote \jpz{}{the} proliferation and death rates of precursor cells and depend on the level of cell maturation. 
Mature cells do not divide and die at the rate $\mu$.  The whole process is regulated by a single feedback mechanism based on the assumption that  there exist signaling molecules (cytokines) which regulate the differentiation or proliferation process. The intensity of the signal depends on the level of mature cells, and is modeled using the dependence
$$
s=s[v(t)]=\frac{1}{1+kv(t)},
$$
which  can be justified using a quasi-steady state approximation of the plausible dynamics of the cytokine molecules, see \cite{Anna1}. This expression reflects the heuristic assumption that signal intensity achieves its maximum under absence of mature cells and decreases asymptotically to zero if level of mature cells increases.

%MD Change: I erased this since a(x) is no more defined
%The parameters $a_w$ and $p_w$ satisfy the conditions,
%$$
 % a_w=a(0),\qquad  p_w=p(0).
%$$
% 
The concentration of signaling molecules $s(v)$ influences the length of the cell cycle (proliferation rate $p$) and/or the fraction of stem cells self-renewal ($a_w$) as well as the rate of cell maturation ($g$).
\\
In this model, differentiation of stem cells takes place during mitosis. The differentiation of progenitor cells occurs independently of proliferation. In other words, cells undergo continuous transformations between divisions. We call this process maturation. In the terms of the model this means that in an infinitesimal time interval $(t,\, t+dt)$, the following events occur to a cell
of maturity $x$,
\begin{enumerate}
\item either the cell matures to level $x+dx$, which happens with probability
$g[x,v(t)]dt$,
\item or the cell divides into $2$ daughters, which happens with probability
$p(x)dt$,
with other events occurring with probabilities of the order $o(dt)$.
\end{enumerate}

%MD Change
If we stick to the discrete model proposed in \cite{Anna1}, we obtain the following relations
 linking proliferation and maturation
 \begin{equation}\label{truedata}
\left\{ \begin{array}{l}
g(x,v) = 2[1-\frac{a(x)}{1+kv(t)}]p(x) , 
\\
a_w=a(0),\qquad  p_w=p(0), \qquad 0 < a_w=a(0) \leq 1,
\end{array} \right.
\end{equation}
which does not necessarily mean that differentiation can only occur by division; see the discussion at the end of Section \ref{DiscCont}.

This leads to a simplification of the boundary condition \eqref{StructModel1:bound0} that becomes
$$u(0,t)=w(t).$$
%------------------------------------------------------------------
\subsection{Discrete {\em versus} Continuous Models}\label{DiscCont}
%------------------------------------------------------------------

\subsubsection*{General Setting}

In this section we consider the relationship between the structured population model \eqref{StructModel1:w}-\eqref{StructModel1:v} and the multicompartmental model introduced in  \cite{Anna1}. Following reference \cite{Anna2}, the multicompartmental model can be formulated in a general way,
\begin{eqnarray}\label{GeneralModel}
\frac{d }{dt} u_1&=&p_1(s)u_1 - g_1(u_1,s) - d_1u_1,\\
\frac{d }{dt} u_i&=&p_i(s)u_i + g_{i-1}(u_{i-1},s) - g_i(u_i,s) -  d_i u_i,\:\ \text{for} \:\: i=1,...,n-1 \\
\frac{d }{dt} u_n&=& g_{n-1}(u_{n-1},s) - d_nu_n, \label{GeneralModel:un}                                                    
\end{eqnarray} 
where $g_i(u_i,s)$ denotes a flux of cells from the subpopulation $i$ differentiating to the subpopulation $i+1$. The terms $p_i(s)u_i $ and  $d_iu_i$ describe  cell fluxes due to proliferation and death, respectively. In the general case proliferation or differentiation may depend on signal intensity. 

%MD Change
In reference \cite{Anna1}, differentiation was linked to proliferation and the following expressions were proposed:
\begin{equation}\label{def:gi}
g_i(u_i,s)=2(1 - \f{a_i}{1+ku_n})p_iu_i.
\end{equation}
This is similar to what was defined for stem cells with $p_w$ and $a_w$; here  $\f{a_i}{1+ku_n}$ represents the fraction of cells remaining at the same stage of differentiation $i$ as the parent cell while the $1-\f{a_i}{1+ku_n}$ fraction of cells differentiates to the higher stage $i+1.$
Formulation (\ref{GeneralModel})-(\ref{GeneralModel:un}) \jpz{allows to describe}{describes} the differentiation \jpz{processes}{process independently} \jpz{as a process independent on}{of} cell proliferation \jpz{, which might take place also between the divisions of the cells.}{in the sense that cells either multiply at stage  $i$ or differentiate from compartment $i$ to $i+1$ and so forth.}
\jpz{ 
Assuming that cell differentiation is faster than the process of cell divisions
 and rescaling terms $g_i$ respectively,  we should in limit obtain a structured population model given in the next paragraph.}{
Assuming that cell differentiation occurs at a properly-chosen time scale compared to the time scale of the cell division process, we show in the next paragraph how to obtain, after a suitable renormalization, the structured population model of the next paragraph.
}

\subsubsection*{Continuous Limit}
\label{subsec:discretetocontinuous}
Let us write System \eqref{GeneralModel}--\eqref{GeneralModel:un} in a dimensionless way. We define  $\P,$  $\D,$ $\G_1,$ $\G,$ $\U_1,$ $\U$ and $\U_n$ as characteristic values for the quantities $p_i,$ $d_i,$ $g_1,$ $g_i$ for $i\geq 2,$ $u_1,$ $u_i$ for $2\leq i\leq n-1,$ and $u_n$, respectively. Then, we define dimensionless quantities by $\bar p_i = \f{p_i}{\P}$ etc. 
We make the following hypothesis $$g_i(u_i,s)=g_i(s)u_i.$$  
Then, system \eqref{GeneralModel}--\eqref{GeneralModel:un}  becomes
\begin{eqnarray}\label{GeneralModel:adim}
\frac{d }{dt} \bar u_1&=&{\P}  \bar p_1(s)\bar u_1 - \G_1 \bar g_1 \bar u_1 - \D  \bar d_1 \bar u_1,\\
\frac{d }{dt} \bar u_2&=&\P \bar p_2(s)\bar u_2 + \G_1 \f{\U_1}{\U} \bar g_1 \bar u_1 - \G \bar g_2 \bar u_2 -  \D \bar d_2 \bar u_2,\label{GeneralModel:u2:adim} \\
\frac{d }{dt} \bar u_i&=&\P \bar p_i(s)\bar u_i + {\G}(\bar g_{i-1}\bar u_{i-1} - \bar g_i\bar u_i) -  \D \bar d_i \bar u_i,\:\:\ \text{for} \:\: i=3,...,n-1, \label{GeneralModel:ui:adim} \\
\frac{d }{dt} \bar u_n&=& \f{\G \U}{\U_n} \bar g_{n-1} \bar u_{n-1} - \D \bar d_n\bar u_n. \label{GeneralModel:un:adim}                                                    
\end{eqnarray} 
\jpz{Supposing that}{Letting} the number of compartments \jpz{tends}{tend} to infinity, we pass from the discrete model  to the continuous model by associating to the $u_i$'s a function, constant  on intervals of type $(\ep i,\ep(i+1))$, with $\ep\to 0,$ $i\to\infty$ and the product $\ep i$ remaining positive and finite, \jpz{denoting}{say} $n=n_\ep,$ with $\ep n \to x^* \in \R_+^*=(0,+\infty).$ 
\jpz{Constants would then be considered as tending to continuous functions}{Compartment dependent constants tend to continuous functions}, sums over the index $i$ are interpreted as Riemann sums tending to integrals while finite differences give rise to derivatives.
\jpz{ (we}{A precise discussion of the limiting process is outside the scope of this presentation. We} refer, for instance, to \cite{GV,DGL} for recent examples of how to obtain such limits  based on moments estimates.

In order to interpret the terms ${\G} (\bar g_{i-1} \bar u_{i-1} - \bar g_i \bar u_i)$ in Equation \eqref{GeneralModel:ui:adim} and $\G_1 \f{\U_1}{\U} \bar g_1 \bar u_1 - \G \bar g_2 \bar u_2$ in Equation \eqref{GeneralModel:u2:adim}, as a finite differences tending to a derivative, we take 
$$
{\G} = \f{1}{\ep}, \qquad \G_1 \f{\U_1}{\U}=\f{1}{\ep}.
$$
Assuming that $\P \bar p_i(s) \bar u_i$ and $\D \bar d_i \bar u_i$ \jpz{could tend}{tend} toward limits $p(x,s)u(x,t)$ and $d(x)u(x,t)$ leads to
$$\P=1,\qquad \D=1.$$ 
In order to obtain that $\f{\G\U}{\U_n} \bar g_{n-1} \bar u_{n-1} $ converges to the limit $g(x^*)u(t,x^*)$ in Equation \eqref{GeneralModel:un:adim}, we require
$$
1=\f{\G\U}{\U_n}=\f{1}{\ep}\f{\U}{\U_n}\qquad \U_n=\f{\U}{\ep}.
$$
This means that the order of magnitude of the number of mature cells is much larger than the one of the maturing cells and stem cells. This can be interpreted by the fact that $u_i$ tends to $u(t,x)$ a \emph{density} of cells per unit of maturity, whereas $u_1$ and $u_n$ are \emph{numbers} of cells. This scaling follows also from mass balance considerations: System \eqref{GeneralModel}--\eqref{GeneralModel:un} leads to the following mass balance
\begin{equation}\label{eq:MassBalance}
 \f{d}{dt} u_1 + \f{d}{dt} \sum\limits_{i=2}^{n-1} u_i + \f{d}{dt} u_n = p_1 u_1 + \sum\limits_{i=2}^{n-1} p_i u_i - d_1 u_1 - \sum \limits_{i=2}^{n-1} d_i u_i - d_n u_n, 
\end{equation}
meaning that the exchange among compartments at maturation rate $g_i$ does not influence the total growth of the  population. In this equation, we keep the specific values $u_1$ and $u_n$ and interpret the sum for $2\leq i\leq n-1$ as an integral. We obtain
$$\f{d}{dt} u_1 + \f{d}{dt} \sum\limits_{i=2}^{n-1} u_i + \f{d}{dt} u_n=\f{d}{dt} \U_1 \bar u_1 + \U \f{d}{dt} \sum\limits_{i=2}^{n-1}  \bar u_i + \U_n \f{d}{dt}  \bar u_n,$$
which leads to the following choice
$$\f{\U}{\U_n}=\f{\U}{\U_1}=\ep.$$
With this choice and the previous relations, we are led to choose $\G_1=1,$ which allows a limit for Equation \eqref{GeneralModel:adim}. We note however that this implies a different order of magnitude for $\G_1$ and for $\G;$ the interpretation could be that we have divided the previous discrete compartments into smaller ones, of size $\ep,$ where division does not occur but where maturation occurs. In this framework, $\G_1$ is not homogeneous to $\G$ but rather to the integral of $\G$ over a small compartment of size $\ep.$ 

Under these assumptions, let us set
 \[\chi_i^\ep(x)=\chi_{[i\ep,(i+1)\ep)}(x),\]
  with $\chi_A$ being the indicator function of a set $A$. 
  We introduce the piecewise constant function   
\[u^\ep(x,t):=\dis\sum_{i=1}^{n_\ep}  u_i(t)\chi_i^\ep(x).\]
By the same token, we associate  the following functions to the coefficients
\[\begin{array}{l}
d^\ep(x):=\dis\sum_{i=1}^{n_\ep} d_i\chi_i^\ep(x),
\qquad
p^\ep(x,s):=\dis\sum_{i=1}^{n_\ep}  p_i(s)\chi_i^\ep(x),
\qquad
g^\ep(x,u(x),s):=\dis\sum_{i=1}^{n_\ep}  g_i(u_i,s) \chi_i^\ep(x).\end{array}\]
We make the following continuity assumptions on the \jpz{adimensioned}{dimensionless} system:
\begin{equation}\label{as:scaling}\begin{array}{l}
\exists K>0\quad \text{s.t.} \quad |g_i|+|d_i|+|p_i| \leq K,\quad |g_{i+1}-g_i|+|d_{i+1}-d_i| + |p_{i+1}-p_i| \leq \f{K}{i}
\\ \\
p_i, g_i \text{ are uniformly continuous with respect to the variable }s.
 \end{array}
\end{equation}
We define the piecewise constant functions  $g^\ep,$ $d^\ep$ and $p^\ep$ on the respective basis of the discrete coefficients $g_i,$ $d_i$ and $p_i$, similarly  as $u^\ep$ was defined for $u_i.$ Assumption \eqref{as:scaling} leads to their convergence (up to subsequences) to continuous functions $g,$ $d$ and $p$ of both variables $x$ and $s$ (see Lemma 1 of \cite{DGL} for instance).  We can prove (based e.g. on \cite{GV,DGL}) the following result.
\begin{proposition}
\jpz{Let us s}{S}uppose that $u_i^\ep$ is \jpz{}{a} solution of Equation \eqref{GeneralModel:ui:adim} verifying \jpz{$u_i(0) \in l^1.$}{$\left( u_i(t=0)\right) \in l^1.$} Under Assumption \eqref{as:scaling}, for all $T>0,$ there exists a subsequence of $\left(u_i^\ep\right)$ converging towards a limit $u\in{\cal C}(0,T;{\cal M}^1([0,x^*])-\text{weak}-*)$ solution of Equation \eqref{StructModel1:u}, $u_1^\ep$ to a limit $u_1 \in {\cal C}(0,T)$ solution of Equation \eqref{StructModel1:w} where $s=s(v)$ with $v$ a limit of a subsequence of $u_n^\ep.$  The boundary condition \eqref{StructModel1:bound0} is satisfied in a distributional sense and  Equation~\eqref{massbalance1} is satisfied, what is equivalent to a weak formulation of the boundary condition~\eqref{StructModel1:v}.
\end{proposition}

\begin{figure}[htbp]
\begin{center}
\begin{minipage}{17cm}
\includegraphics[width=8cm, height=7cm]{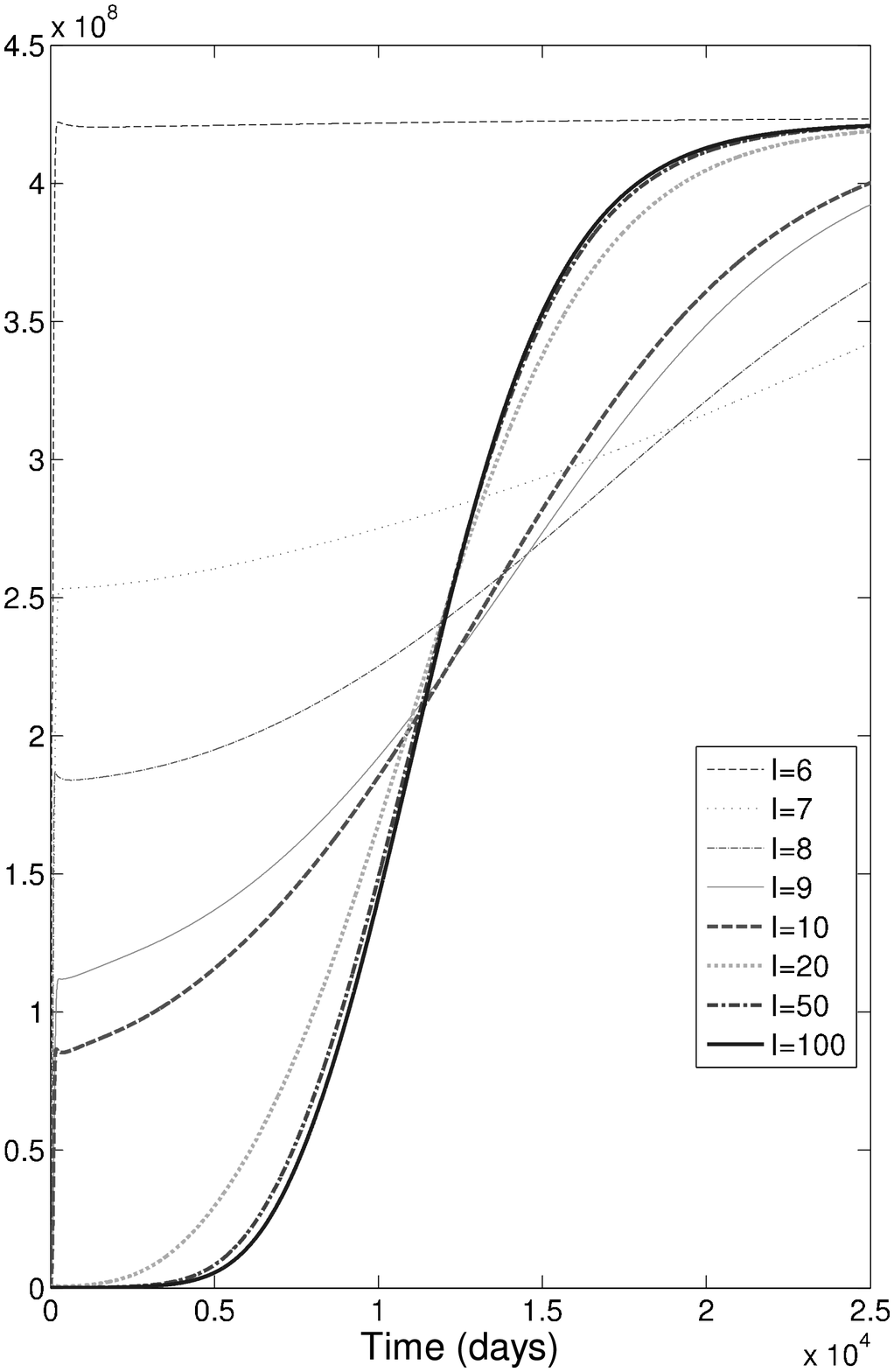} \quad \includegraphics[width=8cm,height=7cm]{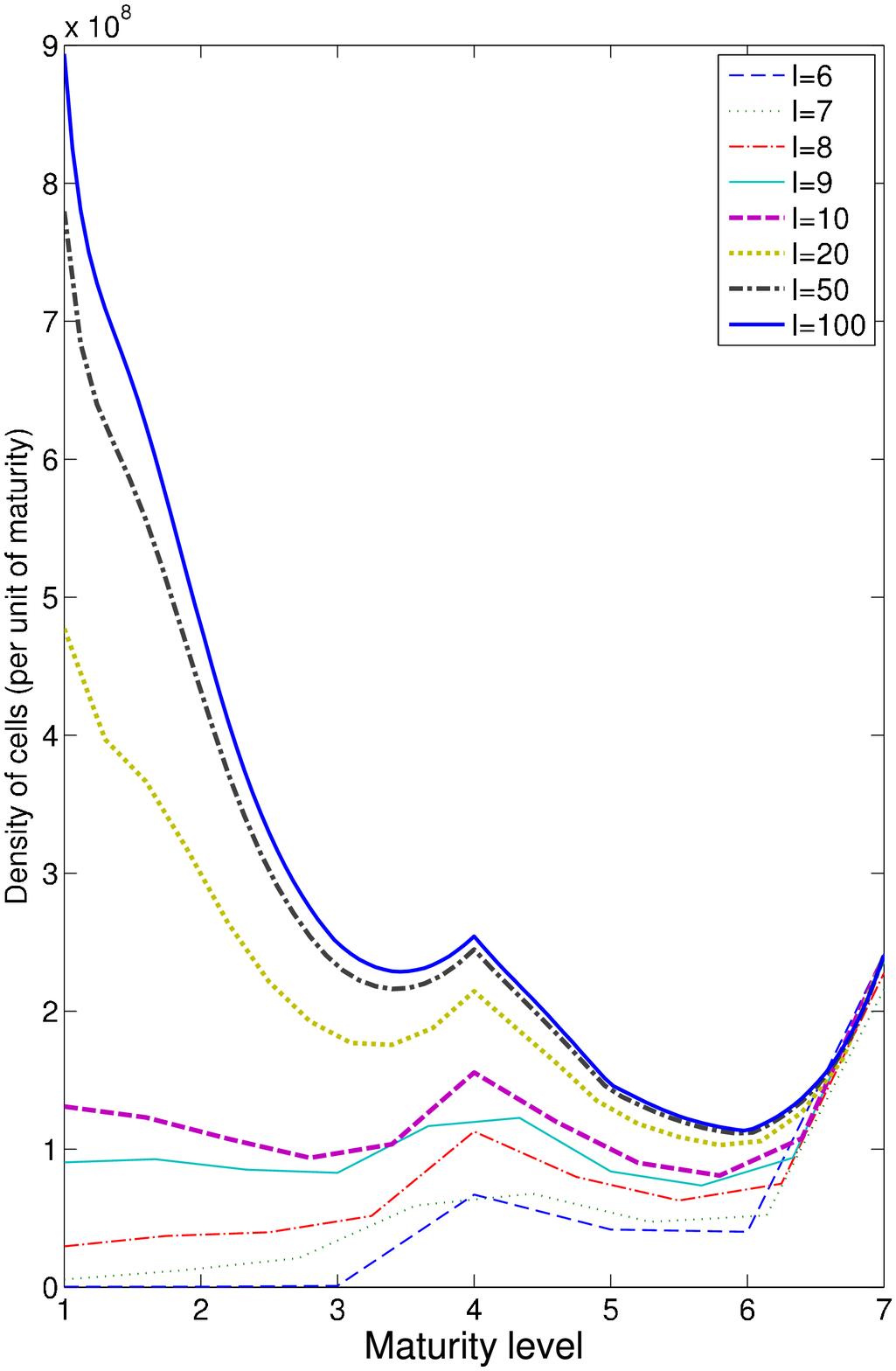}\end{minipage} \end{center}\vspace{-0.5cm}
\caption{ \label{fig:discrete} Comparison of numerical simulations using different grids on the interval $[1,7],$ from $I=6$ (7 points, maturity step $dx=1,$ discrete model) to $I=100$. Left: mature cells evolution with time. Right: distribution of cell density along the maturation level, at steady state. One can see that the model is extremely sensitive to the number of steps (even $7$ to $10$): small numbers seem to be unstable, whereas for large numbers the numerical scheme converges.}
\end{figure}

We note that we have used the fact that maturation and proliferation are decorrelated. If this were not the case, it would be impossible to make  a $1/ \ep$ factor appear in    $\G$, since   $\P=1$. In such case in the limit equation the transport appears as a first order corrective term and Equation \eqref{StructModel1:u} is replaced by
\begin{equation}\label{epstransport}
\partial_t u(x,t)+\ep \partial_x[g(x,s)u(x,t)]  = p(x,s)u(x,t) - d(x)u(x,t).
\end{equation}
Figure \ref{fig:discrete} in Section \ref{subsec:simul} depicts related numerical simulations.

%----------------------------------------------------------------------
\section{Uniform Bounds for the Continuous Model}\label{Analysis}
%----------------------------------------------------------------------
%----------------------------------------------------------------------

In the remainder of this work we will consider  a special version of the above model assuming time independent 
proliferation rates $p(x)$,  and zero death rates of undifferentiated cells $d_w=0$ and $d(x)=0$. Indeed, neglecting death rates of immature cells does not change the analysis. Concerning the feedback loops it was shown in \cite{Anna1} for the discrete model that the feedback on the stem cells self-renewal fraction and on the maturation speed $g$ is much more important for the efficiency of the process than the feedback on the proliferation rate $p(x)$. Therefore, in the reminder of this work we focus on the model with regulated self-renewal and maturation.
We also introduce simpler notation which makes it easier for analysis. This yields the following system of differential equations for $t>0$, $x>0$.
\begin{eqnarray}
 \frac{d}{dt} w(t) & = & { \alpha} (v(t)) w(t), \label{eq:w}\\
\partial_t u(x,t)+\partial_x[g(x,v(t))u(x,t)] & = &p(x)u(x,t), \label{eq:u}\\
u(0,t) & = &  w(t),  \label{eq:uBC}\\
\frac{d}{dt}  v(t) & = & g(x^{*},v(t)) u(x^{*},t)-\mu v(t), \label{eq:v}
\end{eqnarray}
together with initial data
\begin{equation} \label{eq:id}
w(0)=w_0\geq 0,\qquad  u(0,x)=u_{0}(x)\geq 0,\qquad  v(0)=v_{0}\geq 0.
\end{equation}
We obtain the cell number balance law
\begin{equation}\label{massbalance2}
\frac{d}{dt}\big[w +\int u(x,t)dx + v \big]= [ \alpha(v)+g(0,v)] w +\int p(x) u(x,t) dx - \mu v(t),
\end{equation}
corresponding to the fact that the total population can only change by proliferation or death. Indeed, one can interpret $\alpha+g(0,v)$ as the stem cells proliferation rate, see above sections and Equation \eqref{massbalance1}.
%System \eqref{StructModel1:w}--\eqref{StructModel1:v} in Section \ref{Section_model} is the case when 
%\begin{equation}\label{truedata}
%\left\{ \begin{array}{rl}
%{\alpha}(v)&= [2\frac{a_w }{1+kv(t)}-1] p_w,
%\\
%g(x,v) &= 2[1-\frac{a(x)}{1+kv(t)}]p(x) , 
%\\
%p_w&= p(0), \qquad 0 < a_w=a(0) \leq 1,
%\end{array} \right.
%\end{equation}
%but this specific form is not used in the mathematical proofs. 
\\
In the sequel we will study model \eqref{eq:w}--\eqref{eq:v} under the following assumptions
\begin{equation}\label{Ass1}
g_x, \; g_{xx} \in L^\infty([0,x^*]\times \R^+),\;  { \alpha}(x) \in C([0, \infty)),\; p(x)\in C^1([0,x^*]),
\end{equation}
\begin{equation}\label{Ass2}
{ \alpha}(v)\in [\alpha_\infty, \alpha_0], \quad  { \alpha} \; \text{ is  decreasing }, \quad  {\alpha}(+\infty):=\alpha_\infty <0,
\end{equation}
\begin{equation}\label{Ass3}
0< g_- \leq g(x,v)\leq g_+< \infty,  \; \forall (x,v)\in[0,x^*]\times \R^+,
\end{equation}
First we show that the model solutions are uniformly bounded
%--------------------------------------------
\begin{theorem}\label{Bound-wuv}
Under assumptions \eqref{Ass1}--\eqref{Ass3} and that $u_0(x)\in C^1([0,x^*])$, the solution to 
System~\eqref{eq:w}--\eqref{eq:id} is uniformly bounded. More precisely, all the components $w(t)$, $u(x,t)$, $v(t)$ 
are uniformly bounded. 
\end{theorem}
%--------------------------------------------
The remainder of the section is devoted to the proof of this result, which uses some technical lemmas.  
We first prove the following estimate
%-----------------------------------------
\begin{lemma}\label{Bound1}
Under the assumptions of Theorem~\ref{Bound-wuv}, the function   $z(x,t)=\partial_x(\ln {u})$ is uniformly bounded on $[0,x^*]\times \R^+$.
\end{lemma}
%--------------------------------------------
{\bf Proof.}  
The equation for $z$ reads 
\begin{equation}\label{z-equation1st}
\left\{\begin{array}{rl}
\partial_t z + \partial_x(g z)&=-g_{xx}+p_x, \\
z(0,t)&=-\frac{{\alpha}(v)-p(0)}{g(0,v)}-\frac{g_x(0,v)}{g(0,v)} \in L^\infty(0,+\infty).
\end{array}\right.
\end{equation}
Indeed, we have $z(0,t)=  \frac{u_x(0,t)}{u(0,t)}$ and thus we can compute
\begin{eqnarray*}
z(0,t)&=&-\frac{\partial_t u(0,t)-p(0)u(0,t)+g_x(0,v) u(0,t)}{g(0,v)u(0,t)}
\\
&=&-\frac{{\alpha}(v)-p(0)}{g(0,v)}-\frac{g_x(0,v)}{g(0,v)}.
\end{eqnarray*}
And we conclude that $ z(0,t)$ is uniformly bounded by assumptions \eqref{Ass1}--\eqref{Ass3}.\\
Next, we rewrite the equation for $z$ as
\begin{equation}\label{z-equation}
\partial_t z + g \partial_x z = -g_x z + Q(x,v),
\end{equation}
where $Q=-g_{xx}+p_x$ is a bounded function of $v$ and $x$.\\ \\
In the following, we show that the solution to \eqref{z-equation} satisfies the estimate
\begin{equation}\label{z-bound}
||z(x,t)||_{L^\infty} \leq M:= \left( \sup_t| z(0,t) | +\sup_x| z(x,0) |  + x^* \; ||\f{Q}{g}||_{L^\infty} \right)e^{x^* ||\f{g_x}{g}||_{L^\infty}}.
\end{equation}
Indeed, since $g\geq g_- >0,$ we can rewrite Equation~\eqref{z-equation} as
$$ 
\partial_x z +\f{1}{g} \partial_t z  = -\f{g_x}{g} z + \f{Q(x,v)}{g},
$$
and apply the method of characteristics by defining as usual (except that $x$ plays the role of time),
$$
\f{dT}{dx} (x,t) = \f{1}{g}(x,v(T(x,t))),\quad T(x=0,t)=t,
$$
$$ \f{dT^{-1}}{dx} (x,t') = - \f{1}{g}(x,v(T^{-1}(x,t'))),\quad T^{-1}(x=0,t')=t'.$$
We look for solutions of the form
$Z(x,t)=z(x,T(x,t)),$ which satisfy the following equation
$$ \partial_x Z   +\f{g_x}{g} (x,T(x,t)) Z = \partial_x (Z e^{\int\limits_0^x \f{g_x}{g} (\zeta,T(\zeta,t))d\zeta})e^{-\int\limits_0^x \f{g_x}{g} (\zeta,T(\zeta,t))d\zeta}=\f{Q}{g}(x,v(T(x,t))).$$
Integrating the above equation yields
\begin{eqnarray*}
Z(x,t)e^{\int\limits_0^x \f{g_x}{g} (\zeta,T(\zeta,t))d\zeta} &=& Z(0,t) + \int\limits_0^x \f{Q}{g}(\xi,v(T(\xi,t))) e^{\int\limits_0^\xi \f{g_x}{g} (\zeta,T(\zeta,t))d\zeta} d\xi,\\
Z(x,t) &=& Z(0,t)e^{-\int\limits_0^x \f{g_x}{g} (\zeta,T(\zeta,t))d\zeta} + \int\limits_0^x \f{Q}{g}(\xi,v(T(\xi,t))) e^{-\int\limits_\xi^x \f{g_x}{g} (\zeta,T(\zeta,t))d\zeta}  d\xi,\\
z(x,T(x,t)) &=& z(0,t)e^{-\int\limits_0^x \f{g_x}{g} (\zeta,T(\zeta,t))d\zeta} + \int\limits_0^x \f{Q}{g}(\xi,v(T(\xi,t))) e^{-\int\limits_\xi^x \f{g_x}{g} (\zeta,T(\zeta,t))d\zeta}  d\xi.
\end{eqnarray*}
Defining $\bar t= T(x,t),$ or yet $t=T^{-1} (x,\bar t),$  yields $T^{-1}\geq 0$ for $\bar t \geq \f{x^*}{g_{min}},$ and we obtain
$$z(x,\bar t) = z(0,T^{-1} (x,\bar t))e^{-\int\limits_0^x \f{g_x}{g} (\zeta,T(\zeta,T^{-1} (x,\bar t)))d\zeta} + \int\limits_0^x \f{Q}{g}(\xi,v(T(\xi,T^{-1} (x,\bar t)))) e^{-\int\limits_\xi^x \f{g_x}{g} (\zeta,T(\zeta,T^{-1} (x,\bar t)))d\zeta}  d\xi.$$
Therefore, for $\bar t \geq \f{x^*}{g_{min}}$ it holds
$$||z||_{L^\infty} \leq (|z(0,\cdot)|  + x^*||\f{Q}{g}||_{L^\infty} )e^{x^* ||\f{g_x}{g}||_{L^\infty} }.$$
\qed
\\
From Lemma~\ref{Bound1}, we deduce several useful estimates  
%---------------------------------------------
\begin{lemma}\label{Bound2}
There exist positive constants $M_1$, $M_2$, $M_3$ such that the solutions to system \eqref{eq:w}--\eqref{eq:id}  satisfy
\begin{itemize}
 \item[(i)]    $w(t)\leq M_1u(x,t)$,
 \item[(ii)]   $w(t)\leq M_2 v(t)$,
 \item[(iii)]  $u(x,t)\leq M_3 w(t)$.
\end{itemize}
\end{lemma}
%---------------------------------------------

\noindent{\bf Proof} 
$(i)$ Boundedness of $-z=-\frac{\p}{\p x} \ln u$ results in the following inequality
$$
\ln \frac{1}{u} \leq \ln \frac{1}{w} +Mx,
$$
which in turn yields assertion $(i)$ with $M_1=e^{Mx^*}$.
\\ \\
$(ii)$ To bound $w$ by $v,$ we calculate
$$
\f{d}{dt} \f{w}{v}=\f{w}{v} \left({ \alpha}(v(t))  - g(x^*,v(t)) \f{u(x^*,t)}{v} +\mu \right).
$$
Since ${\alpha}(v) \leq { \alpha}(0)$, $g(x^*,v(t))\geq g_-$ and $u(x^*,t)\geq w(t)/M_1$ we obtain
$$
\f{d}{dt} \f{w}{v}  \leq \f{w}{v} \left( {\alpha}(0) +\mu - \frac{g_- }{M_1} \f{w}{v} \right).
$$
This yields the estimate
$$
w(t) \leq v(t) \max \left(  \f{w(0)}{v(0)}, M_1 \frac{{\alpha}(0) +\mu}{ g_-} \right):=M_2 v(t),
$$
and the assertion $(ii)$ is proved.
\\ \\
$(iii)$ The proof follows as in $(i)$, departing from $ \ln u (x,t)  \leq \ln {w(x,t)} +Mx$. 
\qed

As a consequence of Lemma \ref{Bound2}, we derive
%---------------------------------------------
\begin{corollary} 
Under the assumptions of Theorem~\ref{Bound-wuv}, the components $w(t)$, $u(x,t)$   and $v(t)$ of the solutions to 
System~\eqref{eq:w}--\eqref{eq:id} are uniformly bounded.
\end{corollary}
%---------------------------------------------
{\bf Proof.}  Applying Lemma \ref{Bound2} $(ii)$  to equation \eqref{eq:w}, we obtain
 $$
 \frac{dw}{dt}\leq  { \alpha} \left( \frac{ w}{M_2}\right)  w.
 $$ 
This yields boundedness of $w$ by Assumption~\eqref{Ass2}.

 Boundedness of $w$ yields also  boundedness of $u$ using  Lemma \ref{Bound2} $(iii)$. Finally, boundedness of $v$ results from Equation \eqref{eq:v} due to boundedness of $u(x^*,t)$ because $g\leq g_+$. 
\qed
\\
The proof of Theorem \ref{Bound-wuv} is now complete. \qed
 \\
We also state another result, in the spirit of Lemma~\ref{Bound2}, that is used later on
%---------------------------------------------
\begin{lemma}\label{Bound3} 
There exists a constant $M_4>0$  and $0<\gamma<1$ such that $v(t)\leq M_4 w^{\gamma}(t)$. 
\end{lemma}
\noindent {\bf Proof.}
We calculate 
\begin{eqnarray*}
\f{d}{dt} \f{v}{w^\gamma} \leq M_3 g(x^*,v(t))w^{1-\gamma} - \f{v}{w^\gamma}(\mu+\gamma{ \alpha}(v) ).
\end{eqnarray*}
We choose  $\gamma >0$ small enough such that $\mu  + \gamma \alpha_\infty:= \mu_1 >0$ and $\gamma <1$. Since $w$ is uniformly bounded, we find
\begin{eqnarray*}
\f{d}{dt} \f{v}{w^\gamma} \leq C - \f{v}{w^\gamma}\mu_1
\end{eqnarray*}
which yields boundedness of $\f{v}{w^\gamma}$. 
\qed
\\
Finally, we conclude this section with a consequence of Theorem~\ref{Bound-wuv}.
%-----------------------------------------------------
\begin{corollary}\label{Existence}
Under the Assumptions~\eqref{Ass1}--\eqref{Ass3} and that $u_0(x)\in C^1([0,x^*])$, System~\eqref{eq:w}--\eqref{eq:id}  has a unique global solution. 
Furthermore, such solution is uniformly bounded.
\end{corollary}

\noindent {\bf Proof.} Local in time existence of the unique solution follows from the Cauchy-Lipschitz theorem. Theorem \ref{Bound-wuv} provides uniform boundedness of solutions and hence the global existence. 
\qed

%-------------------------------------------

%-------------------------------------------------------
\section{Extinction and Persistence }\label{Section_extinct}
%-------------------------------------------------------

In this section we provide conditions for extinction and persistence of positive solutions.
\\

First, we consider a case when ${ \alpha}(0) < 0$.  In this case there exists only a trivial steady state of the model and
%--------------------------------------
\begin{theorem}
Assume \eqref{Ass1}--\eqref{Ass3}. If ${\alpha}(0) < 0$, then all solutions of  system \eqref{eq:w}--\eqref{eq:id} converge to zero at an exponential rate.
\end{theorem}
%---------------------------------------
 
\noindent{\bf Proof}
First of all, notice that, since ${\alpha}(v) \leq { \alpha}(0)<0$, it is obvious 
from equation \eqref{eq:w}  that $w$ converges to $0$ exponentially.
\\
For the other components, we consider a functional 
$\gamma w(t) +\int_0^{x^{\ast}}e^{-\beta x} u(x,t)dx +  e^{-\beta x^{\ast} }v$, with positive 
constants $\gamma$ and $\beta$ to be determined. We compute its time derivative,
\begin{eqnarray*}
&& \frac{d}{dt}\big( \gamma w(t) +\int_0^{x^{\ast}}e^{-\beta x} u(x,t)dx +  e^{-\beta x^{\ast} }v\big)\nonumber \\
&=&\gamma { \alpha}(v)  w(t) - \beta \int_0^{x^{\ast}}  e^{-\beta x} g(x,v) u(x,t) dx  - e^{-\beta x^{\ast}} g(x^{\ast}, v)u(x^{\ast},t) \\ 
   &&  \; +  g(0,v) u(0,t) + \int_0^{x^{\ast}} p(x) u(x,t) e^{-\beta x}  dx  +  
g(x^{\ast}, v) u(x^{\ast},t) e^{-\beta x^{\ast} }- \mu e^{-\beta x^{\ast} }v \nonumber \\
&=&[\gamma { \alpha}(v) +  g(0,v) ] w(t) +  \int_0^{x^{\ast}}  e^{-\beta x}  u(x,t)\big(p(x)- \beta g(x,v)\big )dx - \mu e^{-\beta x^{\ast} }v.
\end{eqnarray*}
Since ${ \alpha}(v) \leq{ \alpha}(0)<0$ we may choose $\gamma$ such that 
$\gamma { \alpha}(0)  + \sup_vg(0,v) \leq - \Gamma <0$. Moreover, choosing $\beta$ such that $p(x)- \beta g(x,v)<-\Gamma$ implies that
\begin{eqnarray*}
\frac{d}{dt}\big( \gamma w(t) +\int_0^{x^{\ast}}e^{-\beta x} u(x,t)dx +  e^{-\beta x^{\ast} }v\big)\leq -\Gamma  w(t) - \Gamma \int_0^{x^{\ast}}e^{-\beta x} u(x,t)dx- \mu e^{-\beta x^{\ast} }v.
\end{eqnarray*}
We conclude that the solutions converge to zero at an exponential rate for  $t \rightarrow \infty$.
\qed

Secondly, if it is the case that ${ \alpha}(0)>0$, then we conclude that the solutions to the system 
cannot become extinct.
%--------------------------------------------
\begin{theorem}\label{nonext}
Assume \eqref{Ass1}--\eqref{Ass3}, $w(0)>0$ and $u_0(x)\in C^1([0,x^*])$. If  ${ \alpha}(0)>0$,  the solution $u$, $v$, $w$ of system \eqref{eq:w}-\eqref{eq:v}  with positive initial conditions remain bounded away from zero.
\end{theorem}
%--------------------------------------------

\noindent {\bf Proof}
Applying Lemma \ref{Bound3} to equation \eqref{eq:w}, we obtain
$$
 \frac{dw}{dt}\geq  {\alpha} \left( M_4 w^\gamma \right) w,
$$ 
and the assumption ${ \alpha}(0)>0$ allows us to conclude. Then, the estimates of Lemma~\ref{Bound2} conclude for $u$ and $v$.
\qed

%----------------------------------------------------------------------------
\section{Stationary Solutions and Their Stability}\label{Section_stationary}
%-------------------------------------------

\subsection{Stationary Solutions}
As usual in dynamical systems, a natural question concerns the existence of steady states (stationary solutions). We shall
now investigate this issue.

In our case, the steady states are given by the solutions $(\bar{w}, \bar{u}, \bar{v})$  to the system
\begin{eqnarray}
 { \alpha}(\bar v) \bar w& =&0,\label{SteadyState1}\\
\frac{d}{dx}[\bar g(x) \bar u(x)] & = &p(x)\bar u(x),\label{SteadyState2}
\\
\bar u(0) &=& \bar w, 
\\
  \bar g(x^{*}) \bar u(x^{*})-\mu \bar v& = &0,\label{SteadyState3}   
\end{eqnarray} 
where $\bar g(x) := g(x,\bar v)$. 

System~\eqref{eq:w}--\eqref{eq:v} always admits the trivial steady state $w=0$, $u=0$, $v=0$, 
which we do not consider.
Depending upon the value ${ \alpha}(0)$ it  may also have exactly one positive steady state $(\bar w, \bar u,\bar v)$, as we state it in the 

%-------------------------------------------
\begin{lemma}\label{lm:stst}
Under the Assumptions~\eqref{Ass1}--\eqref{Ass3}, the System~\eqref{eq:w}-\eqref{eq:v} has  
a strictly positive steady state if and only if  ${\alpha}(0)>0$.  Furthermore, the steady state is unique.
\end{lemma}
%-------------------------------------------

This condition is in agreement with biological observations concerning self-renewal of stem cell subpopulation \cite{He} and an analogous condition for the compartmental model was discussed in \cite{Anna2}.
\\ \\
{\bf Proof} 
Since we discard the trivial steady state, from equation (\ref{SteadyState1}) we obtain the condition  ${ \alpha}(\bar v)=0$. 
As we know that ${ \alpha}$ decreases and tends to $\alpha_\infty<0$ at infinity, there exists a unique  solution $\bar v$ to 
\begin{equation}\label{StationaryV}
{ \alpha}(\bar v)=0,
\end{equation}
 if and only if the condition ${ \alpha}(0)>0$ holds. 
Thus, we may compute $$\bar u(x^{*}) = \f{\mu \bar v}{\bar g(x^*)} \mbox{ .}$$
% \\
Solving differential equation (\ref{SteadyState2}) with this boundary condition at $x=x^*$ yields,
\begin{equation}\label{StationaryU}
\bar u(x)= \frac{\bar g(x^*)}{\bar g(x)} \bar u(x^*) \exp\Big\{-\int_{x}^{x^*} \frac{p(\xi)}{\bar g(\xi)}d\xi\Big\}.
\end{equation}
We finally identify $\bar w$ using the boundary condition at $x=0$, $\bar w = \bar u(0)$, which leads to 
\begin{equation}\label{StationaryW}
\bar w = \frac{\bar g(x^*)}{\bar g(0)} \bar u(x^*) \exp\Big\{-\int_{0}^{x^*} \frac{p(\xi)}{\bar g(\xi)}d\xi\Big\} = \frac{\mu \bar v}{\bar g(0)}\exp\Big\{-\int_{0}^{x^*} \frac{p(\xi)}{\bar g(\xi)}d\xi\Big\}.
\end{equation}
This gives explicit values of the model and completes the proof.
\qed
\\
For ${\alpha}$ given explicitly by  \eqref{truedata},  we may compute
\begin{equation}\label{StationaryPart}
\left\{ \begin{array}{rl}
\bar v&=\frac{2a_w-1}{k},
\\[2mm]
\bar u(x^{*}) &= \frac{\mu }{kp(x^{*})} \; \; \frac{a_w(2a_w-1)}{2a_w-a(x^{*})}.
\end{array}\right.
\end{equation}

\subsection{The Linearized Problem around the Steady State}\label{Stability}
%---------------------------------------------------------------------

In order to investigate local linear stability, we consider in this section the linearization around the positive steady state.
We first derive a characteristic equation for the eigenvalue problem. The signs of the real parts of these eigenvalues give stability (if they all are negative) or instability (if there exists one with positive real part). To emphasize our main point, which is that stability as well as instability of the positive steady state can take place for a suitable choice of model parameters, 
we shall focus on some simpler cases where stability analysis is more  transparent.

%-----------------------------------------------------------------------------------
\subsubsection*{The Characteristic Equation in the General Case}
%-----------------------------------------------------------------------------------
  
We denote by $(\bar w, \bar u,\bar v)$ the steady state solution to Equations~\eqref{SteadyState1}--\eqref{SteadyState3}. Positivity of the considered steady state yields
$\alpha(\bar v)=0$. The linearized problem reads
\begin{eqnarray}\label{eq:steady:lin1}
 \frac{d}{dt} w(t)  & = &    \f{d  \alpha}{dv} (\bar v) \bar w v(t),\\
\partial_t u(x,t)+\partial_x[g(x,\bar v)u(x,t)] + \partial_x[\f{\partial g}{\partial v} (x,\bar v) \bar u(x)] v(t)  & = &p(x)u(x,t),\\
u(0,t)&=&   w(t), \\
\frac{d}{dt}  v(t)  =  g(x^{*},\bar v) u(x^{*},t)+ \f{\partial g}{\partial v} (x^{*},\bar v) \bar u(x^{*}) v(t)&-&\mu v(t), \label{eq:steady:lin2}
\end{eqnarray}
where $w$, $u$, and $v$ denote now the deviation of the solution from the steady state.
Setting  $w(t)=We^{\lambda t}$, $u(x,t)=U(x)e^{\lambda t}$ and $v(t)=Ve^{\lambda t}$ we obtain the eigenvalue problem of the form
\begin{eqnarray}\label{eq:steady:eigen1}
 \lambda W  &= & \f{d \alpha}{dv} (\bar v) \bar w V,\\
\lambda U(x)+\partial_x[g(x,\bar v)U(x)] + \partial_x[\f{\partial g}{\partial v} (x,\bar v) \bar u(x)] V  & = &p(x)U(x),\label{eq:steady:eigenU}\\
U(0)&= & W , \label{eq:steady:eigen3}\\
\lambda V  =  g(x^{*},\bar v) U(x^{*})+ \f{\partial g}{\partial v} (x^{*},\bar v) \bar u(x^{*}) V&-&\mu V. \label{eq:steady:eigen2}
\end{eqnarray}
%
%$$(\lambda -p(x)) U(x)+\partial_x[g(x,\bar v)U(x)] = \p_x [g(x,\bar v)U(x) e^{\int\limits_0^x \f{\lambda -p(s)}{g(s,\bar v)}ds}] e^{-\int\limits_0^x \f{\lambda -p(s)}{g(s,\bar v)}ds} =- \partial_x[\f{\partial g}{\partial v} (x,\bar v) \bar u(x)] V,\\$$
Defining an auxiliary function $f(x)$ such that
$$f(x)e^{\int\limits_0^x \f{-p(s)}{g(s,\bar v)} ds}=- \partial_x[\f{\partial g}{\partial v} (x,\bar v) \bar u(x)],$$  
we obtain
$$
\p_x [g(x,\bar v)U(x) e^{\int\limits_0^x \f{\lambda -p(s)}{g(s,\bar v)}ds}] = f(x) Ve^{\int\limits_0^x \f{\lambda }{g(s,\bar v)}ds},
$$
$$
g(x,\bar v) U(x) = g(0,\bar v) U(0)e^{-\int\limits_0^x \f{\lambda -p(s)}{g(s,\bar v)}ds} + V e^{-\int\limits_0^x \f{\lambda -p(s)}{g(s,\bar v)}ds} \int\limits_0^x f(s)e^{\int\limits_0^s \f{\lambda }{g(\sigma,\bar v)}d\sigma}ds.
$$ 
Hence, using \eqref{eq:steady:eigen1} and \eqref{eq:steady:eigen3} leads to
$$
g(x^*,\bar v) U(x^*) = \biggl(g(0,\bar v) \f{d \alpha}{dv} (\bar v)  \f{\bar w}{\lambda}   +    \int\limits_0^{x^*} f(s)e^{\int\limits_0^s \f{\lambda }{g(\sigma,\bar v)}d\sigma}ds\biggr)V e^{-\int\limits_0^{x^*} \f{\lambda -p(s)}{g(s,\bar v)}ds}. 
$$
We insert this expression in Equation \eqref{eq:steady:eigen2} and obtain the characteristic equation
\begin{equation}\label{eq:eigenvalue:general}
\lambda + \mu- \f{d g}{d v} (x^{*},\bar v) \bar u(x^{*})   =  \biggl(g(0,\bar v) \f{d  \alpha}{dv} (\bar v)  \f{\bar w}{\lambda}   +    \int\limits_0^{x^*} f(s)e^{\int\limits_0^s \f{\lambda }{g(\sigma,\bar v)}d\sigma}ds\biggr) e^{-\int\limits_0^{x^*} \f{\lambda -p(s)}{g(s,\bar v)}ds}.
\end{equation}

%-----------------------------------------------------------------------------------
\subsubsection*{The Simplest Case: $g$ independent of $v$}
\label{sec:generalised:originalmodel}
%-----------------------------------------------------------------------------------

We first focus on the simplest case when the maturation rate $g(x,v)=g(x)$ does not depend on $v$. In other words, the feedback loop only affects stem cells. 
Although this case is very restrictive compared to the original discrete model, since it does not include
relation \eqref{truedata}, it is an illustrative example of a possible general behavior.  Instability and appearance of the oscillations in this model suggest that regulation of the processes solely by the stem cell level is not enough to stabilize the system. Moreover, regulatory feedback between mature cells and progenitor cells has a stabilizing effect and is essential for efficient regulation of the process.\\
Since we have $f=0$, combining Equation~\eqref{StationaryW} with Equation~\eqref{eq:eigenvalue:general}, we arrive at
\begin{equation}\label{eq:eigenvalue:gconstant}
\lambda^2 + \mu\lambda =  \mu \bar v \f{d  \alpha}{dv} (\bar v)  e^{- \tau \lambda}, \qquad \tau =\int\limits_0^{x^*} \f{1}{g(s)}ds>0.
\end{equation}
The relationship is identical with the characteristic equation of a delay differential system. Indeed, problem~\eqref{eq:w}--\eqref{eq:v} 
can be reformulated as a delay differential system. We obtain the following result.

%----------------------------
\begin{proposition}\label{prop:instab:simplest}
Assume that Equations~\eqref{Ass1}--\eqref{Ass3} hold, $\alpha (0)>0$, and $g$ is independent of $v$. Consider 
the steady state $(\bar u,\bar v,\bar w)$ given in Lemma \ref{lm:stst}. Then, \\
(i) for $1<\tau \bar v \; |\f{d \alpha}{dv} (\bar v) |< \f{\pi}{2}$, the system undergoes a Hopf bifurcation for a 
single value $\mu_0>0$ of the parameter $\mu$. Therefore the steady state can be either locally stable or unstable. \\
(ii) Further bifurcations also occur for 
$\tau \bar v \; |\f{d  \alpha}{dv} (\bar v) |> 2k \pi + \f{\pi}{2}$ and $k \geq 1$ for at least one value $\mu_k>0$.
\end{proposition}
%----------------------------

Because in the special case at hand, the system can be reduced to a delay differential equation, the linearised stability implies the stability of the nonlinear system, which then undergoes a Hopf bifurcation for certain values of the parameters (see \cite{Diekmann,MagalRuan} and the references therein).\\

\begin{figure}[htbp] 
\begin{center}
\begin{minipage}{17cm}
\includegraphics[width=8cm, height=7cm]{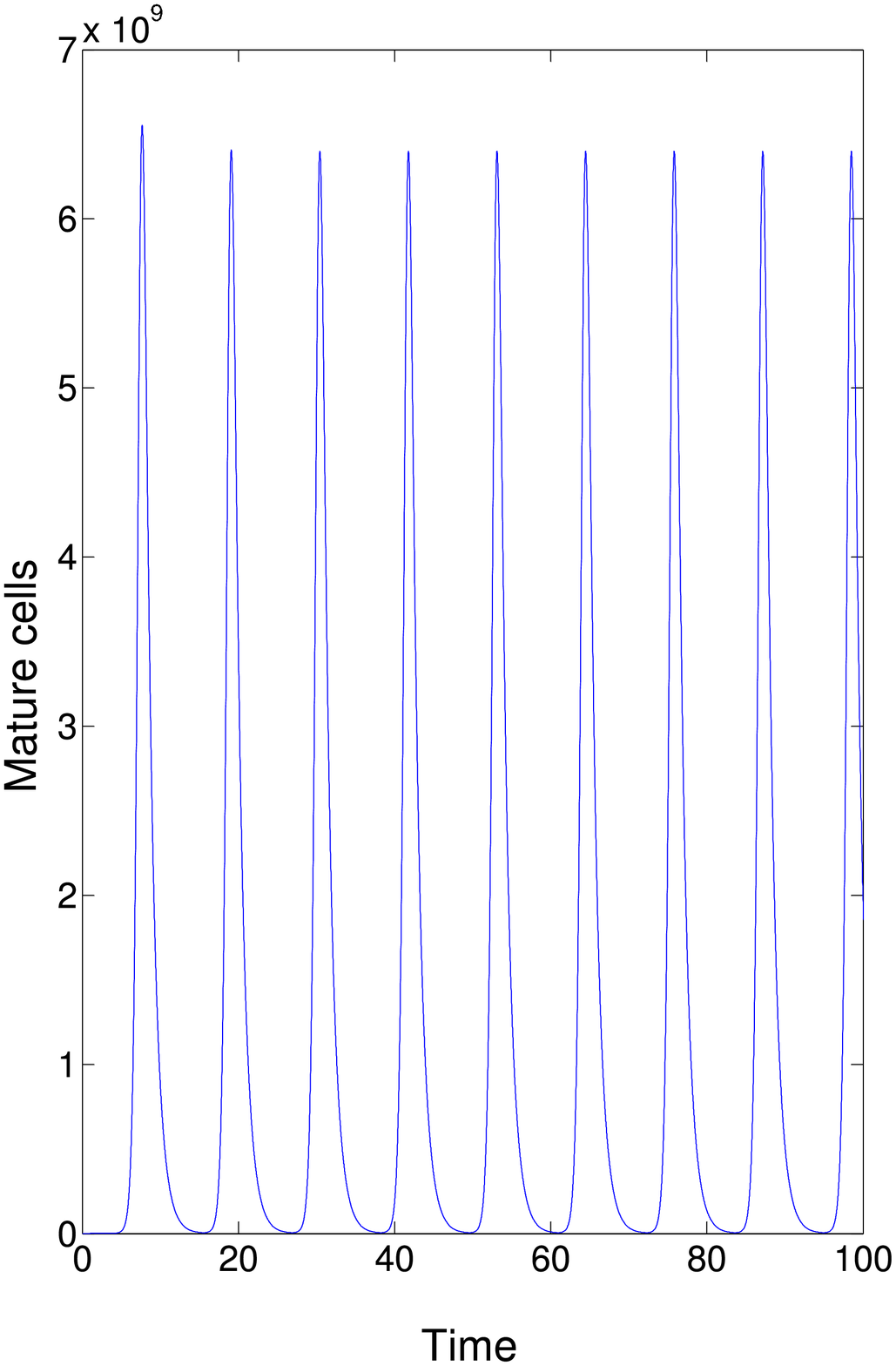} \quad \includegraphics[width=8cm,height=7cm]{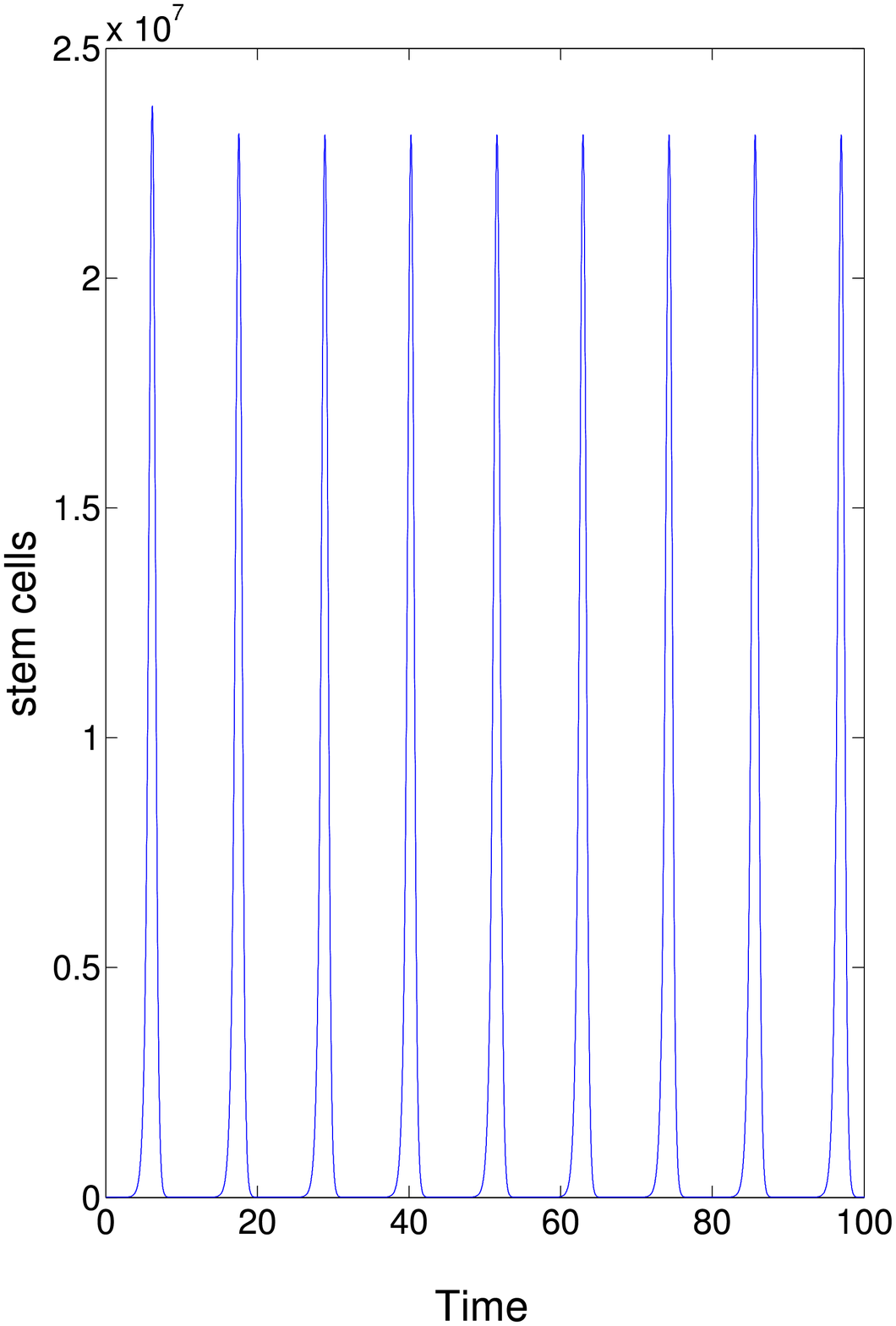}\end{minipage} \end{center}\vspace{-0.5cm}
\caption{\label{fig:instab1:1} Example of instability, in illustration of Proposition \ref{prop:instab:simplest}. Left: evolution of mature cells. Right: evolution of stem cells.}
\end{figure}

\begin{figure}[ht] 
\begin{center}
\begin{minipage}{17cm}
\includegraphics[width=8cm, height=7cm]{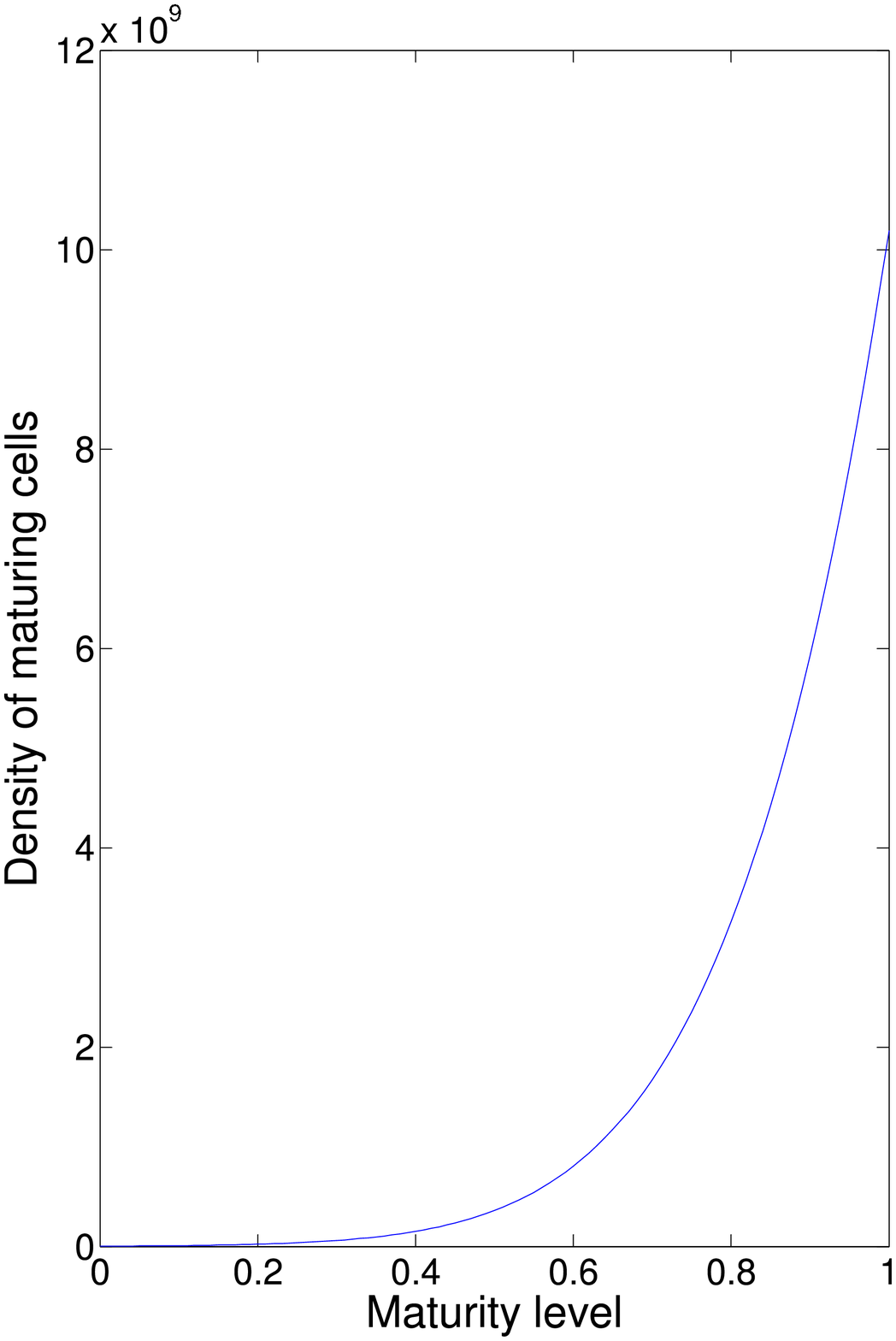} \quad \includegraphics[width=8cm,height=7cm]{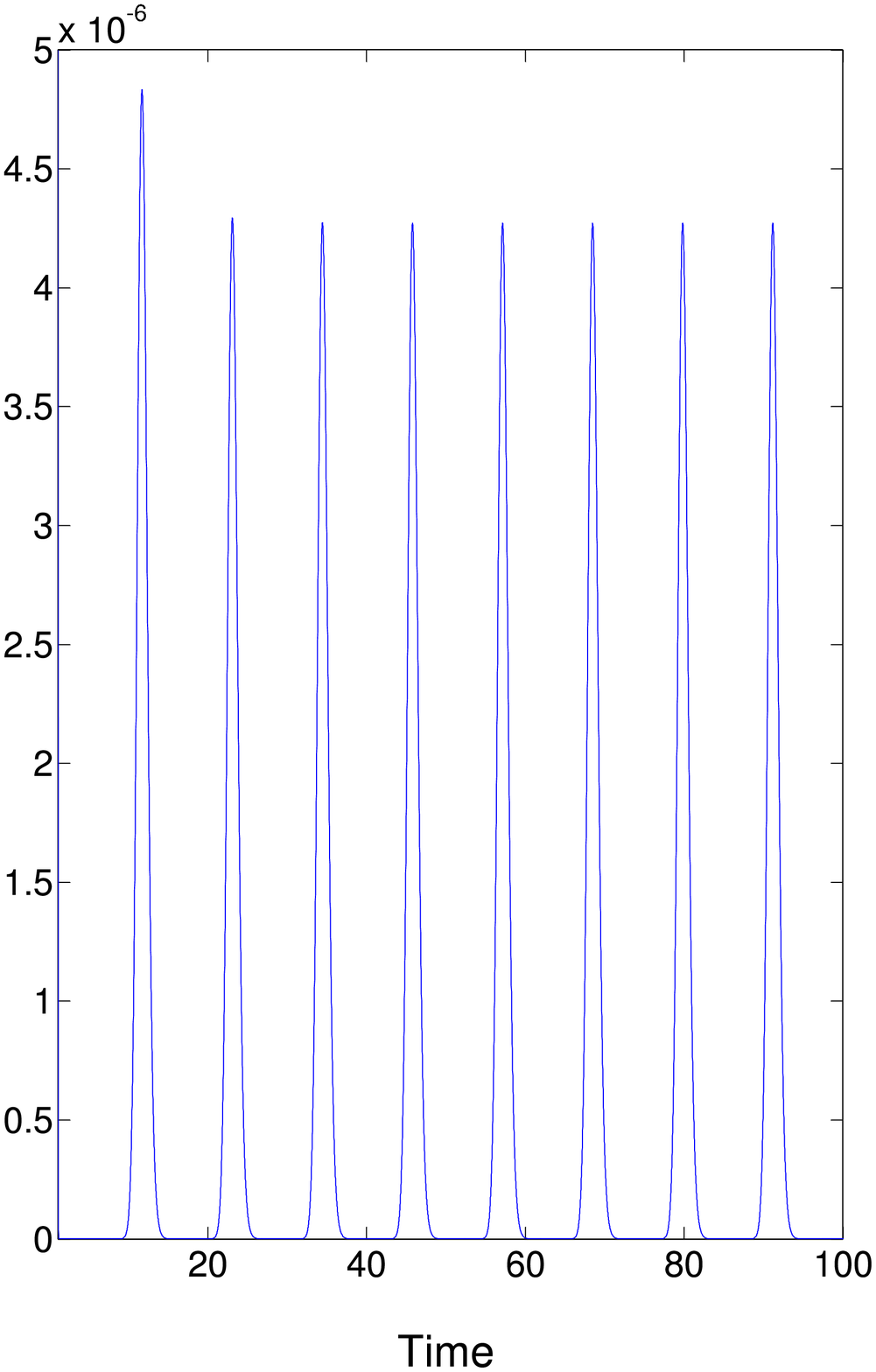}\end{minipage} \end{center}\vspace{-0.5cm}
\caption{\label{fig:instab1:2} Same case as in Figure \ref{fig:instab1:1}. Left: final distribution of cells according to their maturity level. Right: time evolution of $\sqrt{\f{\int |\f{\p}{\p t} u(x,t)|^2 dx}{\int |u(x,t)|^2 dx}}$, to measure the trend to a stable maturity level distribution.}
\end{figure}

\noindent {\bf Proof.}
(i) In order to identify the parameter values for which the bifurcation occurs, 
we look for purely imaginary solutions $\lambda=i\omega$ with $\omega\in \R$.  We obtain the two following relations
$$
\omega^2=\mu  \bar v \; |\f{d  \alpha}{dv} (\bar v) |�\; \cos(\tau \omega),\qquad \tau \omega = \tau \bar v \; |\f{d  \alpha}{dv} (\bar v) |�\; \sin(\tau\omega).
$$
By symmetry, we only consider $\omega >0$. The second relation gives a single value $ \tau \omega_0 \in(0,\frac{\pi}{2})$ as soon as $1<\tau \bar v \; |\f{d  \alpha}{dv} (\bar v) |�< \f{\pi}{2}$. We  can enforce the first relation for a single $\mu$,  because $\cos(\tau \omega_0) >0$. This proves statement (i).
(ii) For $\tau \bar v \; |\f{d  \alpha}{dv} (\bar v) |�> 2k \pi + \f{\pi}{2}$ and $k \geq 1$, the equation $\tau \omega= \tau \bar v \; |\f{d  \alpha}{dv} (\bar v) |� \sin (\tau \omega)$ also has a root $\tau \omega_k \in (2k\pi, 2k \pi + \frac{\pi}{2}$, for which $\cos(\tau \omega_k)>0$ and thus we can find again a $\mu_k$ for which the first equation is satisfied. But there might be multiple compatible crossings and several bifurcations are possible.  
\qed
\\
We now proceed numerically using, for instance Matlab's device DDE BIFTOOL. 
We check that for the values $\mu=\tau=1,$ we get stability for $\mu \bar v \f{d\alpha}{dv}(\bar v)=-1$ and instability for $\mu \bar v \f{d\alpha}{dv}(\bar v)=-2.$
\\
This proposition as well as numerical simulations (see Figures \ref{fig:instab1:1} and \ref{fig:instab1:2}) show that instability occurs through a Hopf bifurcation, and that regular oscillations appear.

%-----------------------------------------------------------------------------------
\subsubsection*{A Case Motivated by the Discrete Model}
\label{sec:linearised:original}
%-----------------------------------------------------------------------------------

In this section we will study more closely the case given by the relations \eqref{truedata}. We can use the values of the steady state computed in Equations \eqref{StationaryV}--%,\eqref{StationaryU*},\eqref{StationaryU} and 
\eqref{StationaryW}, keeping $g(x,v)$ and $p(x)$ fully general. It implies, denoting $\bar g(x)=g(x,\bar v)$
$$\f{d\alpha}{dv} (\bar v) =-\f{2ka_w p_w}{(1+k\bar v)^2}=-\f{kp_w}{2a_w},\qquad \bar w=\frac{\bar g(x^*)}{p_w}\bar u(x^*) exp\Big\{-\int_{0}^{x^*} \frac{p(\xi)}{\bar g(\xi)}d\xi\Big\}.$$
%and so
%\begin{equation}\label{eq:eigenvalue:alpha}
%\lambda + \mu- \mu \f{a(x^*)}{2a_w-a(x^*)}(1-\f{1}{2a_w})   =  \biggl(-\mu \f{p_w}{\lambda} (1-\f{1}{2a_w}) +  e^{\int_{0}^{x^*} \frac{p(\xi)}{\bar g(\xi)}d\xi}  \int\limits_0^{x^*} f(s)e^{\int\limits_0^s \f{\lambda }{g(\sigma,\bar v)}d\sigma}ds\biggr) e^{-\int\limits_0^{x^*} \f{\lambda}{\bar g(s)}ds}.
%\end{equation}
We can show (see the Appendix for detailed calculations) that Equation \eqref{eq:eigenvalue:general} can be written as
\begin{equation}\label{eq:eigenvalue:alpha2}
\lambda + \mu =  \f{\mu}{k}(2a_w-1)\biggl(\f{k}{2a_w}   \bigl(-\f{p_w}{\lambda}+1\bigr)   +    
\int\limits_0^{x^*}  \f{\partial g}{\partial v} (x,\bar v)  \f{\lambda -p(x)}{\bar g(x)^2}e^{\int\limits_0^x \f{\lambda }{\bar g(s)} ds}dx 
\biggr) e^{-\int\limits_0^{x^*} \f{\lambda}{\bar g(s)}ds}.
\end{equation}

\subsubsection*{Case of $\alpha(v)=p_w(\f{2a_w}{1+kv}-1)$ and $g$ independent of $x$}

Since $g(v)$ is now independent of $x$, we have that for all $x$ 
$$g(x,v)=g(0,v)=p_w -\alpha(v)=2p_w(1-\frac{a_w}{1+kv}) \mbox{ .}$$
In particular $g(\bar v)=p_w$ so $\f{dg}{dv} (\bar v)=p_w \f{k}{2a_w}.$ 
We now substitute $g(x,\bar v)=g(\bar v)$ in Equation \eqref{eq:eigenvalue:alpha2} and obtain
$$
\lambda + \mu =  \mu\f{2a_w-1}{2a_w}\biggl(  -\f{p_w}{\lambda}+1   +    
\int\limits_0^{x^*}   \f{\lambda -p(x)}{p_w}e^{\int\limits_0^x \f{\lambda }{p_w} ds}dx 
\biggr) e^{-\int\limits_0^{x^*} \f{\lambda}{p_w}ds}.
$$
We compute the first part of the integral term: $\int\limits_0^{x^*}   \f{\lambda}{p_w}e^{\int\limits_0^x \f{\lambda }{p_w} ds}dx 
=e^{ \f{\lambda }{p_w}x^*}-1,
$
and so
\begin{equation}\label{eq:eigenvalue:alpha3}
\lambda + \f{\mu}{2a_w} =- \f{\mu}{2a_w}(2a_w-1)\biggl(  \f{p_w}{\lambda}   +    
\int\limits_0^{x^*}   \f{p(x)}{p_w}e^{\f{\lambda }{p_w} x}dx 
\biggr) e^{-\f{\lambda}{p_w}x^*}.
\end{equation}
\begin{proposition}\label{prop:instab:originalmodel}
Let $\alpha(v)$ be defined by \eqref{truedata} with $a_w > \f{1}{2},$ and $(\bar u,\bar v,\bar w)$ be defined by 
Equations~\eqref{StationaryV}--\eqref{StationaryW} the unique steady state solution of System \eqref{eq:w}--\eqref{eq:v}. If the maturation rate $g(x,v)$ is independent of the maturity of the cell $x$ and if the proliferation rate $p$ is constant, then the steady state $(\bar u,\bar v,\bar w)$ is locally linearly stable. For a non-decreasing proliferation rate, instability may appear.
\end{proposition}
We treat a case of non-decreasing proliferation rate because it is the most biologically relevant; 
however instability may appear for other cases, 
and is even easier to exhibit, as the proof (postponed to the Appendix) shows. Figures \ref{fig:instab2:1} and \ref{fig:instab2:2} below illustrate a case of instability with a nondecreasing proliferation rate.

\section{Numerical Simulations}\label{Numeric}
In this section we illustrate our theoretical results with a number of numerical simulations. We start with a description of our
numerical methods. 
\subsection{The Numerical Scheme}
We build a simple numerical scheme for System~\eqref{eq:w}--\eqref{eq:v}.
We discretize the problem  on a grid regular in space and adaptive in time. We denote by $\Delta t^k=t^{k+1}-t^k$ the time step between time $t^{k+1}$ and time $t^k,$ by $\Delta x={x^*}/{I}$ the spatial step, 
 where $I$ denotes the number of points: $x_i=i\Delta x,$ $0\leqslant i\leqslant I.$

We use an explicit upwind finite volume method for $u$
$$
u_i^k= \f{1}{\Delta x} \int\limits_{x_{i-\f{1}{2}}}^{x_{i+\f{1}{2}}} u(t^k,y)dy,\qquad  \f{1}{\Delta t^k}\int\limits_0^{\Delta t^k} u(t^k +s, x_{i+\f{1}{2}}) ds  \approx u_i^k.
$$
For time discretization, 
we use a marching technique. At each time $t^k,$
we choose the time step $\Delta t^k$ so as to satisfy the largest possible CFL stability criterion 
$$\theta:= g \f{\Delta t^k}{\Delta x}\leq 1,$$ 
so that
$$\Delta t^k=\f{\Delta x }{Max_x g (x,v^k)}.$$ 
In order to avoid a vanishing time step, it is necessary here to suppose $g \in L^\infty.$ 
Also, more efficient schemes (of WENO type for instance, see \cite{Shu, Jiang}) could be used to capture discontinuities of $g.$

The algorithm is the following: 
\begin{itemize}
\item{\bf Initialization} We use the initial data
$$
w^0=w_0,\qquad  u^0_j=\f{1}{\Delta x} \int\limits_{x_{i-\f{1}{2}}}^{x_{i+\f{1}{2}}} u_0(y)dy,\qquad  v^0=v_{0}.
$$

\item{\bf From $t^k$ to $t^{k+1}:$}

\begin{itemize}
 \item We calculate $\alpha^{k}=\alpha(v^{k})$ and define $w^{k+1}=(1+\Delta t^k \alpha^{k}) w^k.$ 
\item We calculate $\Delta t^k=\f{\Delta x }{Max_i g (x_i,v^k)}$ and define $t^{k+1}=t^k+\Delta t^k.$
\item For a boundary condition at $i=0,$ we define $u_0^{k+1}=w^{k+1}.$
\item We define $u^{k+1}_i$ by the following scheme
$$\frac{u^{k+1}_j - u^k_j}{\Delta t^k} + \frac{g(x_j,v^{k}) u_j^k - g(x_{j-1},v^{k}) u_{j-1}^k}{\Delta x} = p_j u_j^k.$$
\item We define $v^{k+1}$ by
$$\frac{v^{k+1}-v^k}{\Delta t^k} = g(x_I,v^k) u_I^k - \mu v^{k+1},$$
and the term $\mu v^{k+1}$in the right hand side is discretized implicitely for stability. The reason for the choice of $u_I^k$ instead of $u_I^{k+1}$ in the right-hand side of this last scheme
 is due to cell number balance considerations as shown below.
\end{itemize}

\item{\bf Cell number balance.} From Equation \eqref{eq:w}--\eqref{eq:v} we have obtained the cell number balance \eqref{massbalance2}.
 We check the equivalent discrete mass balance:
$$ \f{w^{k+1}-w^k}{\Delta t^k} + \sum\limits_{j=0}^I  \f{u_j^{k+1} -u_j^{k}}{\Delta t^k} + \f{v^{k+1}-v^k}{\Delta t^k}=$$
$$\big(\alpha^k + g(x_0,v^k)\big)  w^k +\sum\limits_{i=0}^I p_j u_j^k \Delta x - \mu v^{k+1}.
$$
\end{itemize}

\subsection{Numerical Simulations}
\label{subsec:simul}

First we compare results of the numerical simulations of the discrete and the continuous models.
To do so, we depart from the discrete values of parameters  given in \cite{ThomasDiploma}.  The notations are those of System \eqref{GeneralModel}-\eqref{GeneralModel:un}, with $g_i(s,u_i) = 2[1-a_i(s)] p_i u_i$, $p_i$ independent of $s,$ $d_i=0$ for $i<n$ and $a_i(s)=\f{a_i}{1+ku_n}$. It corresponds to the model 1 studied in \cite{Anna1, Anna2}. 
\begin{table}[h]\begin{center}
\begin{tabular}{| c| c|| c| c|| c| c|}\hline 
Parameter & Value & Parameter & Value & Parameter & Value  \\ \hline
$a_1$ & $0.77$ & $p_1$ &  $2.15\dot10^{-3}$ day$^{-1}$ & $d_8$ & $0.6925$ day$^{-1}$ \\ \hline
$a_2$ & $0.7689$ & $p_2$ &  $11.21\dot10^{-3}$ day$^{-1}$ & $k$ & $12.8.10^{-10}$  \\ \hline
$a_3$ & $0.7359$ & $p_3$ &  $5.66\dot10^{-2}$ day$^{-1}$ & &  \\ \hline
$a_4$ & $0.7678$ & $p_4$ &  $0.1586$ day$^{-1}$ & &  \\ \hline
$a_5$ & $0.154$ & $p_5$ &  $0.32$ day$^{-1}$ & &  \\ \hline
$a_6$ & $0.11$ & $p_6$ &  $0.7$ day$^{-1}$ & &  \\ \hline
$a_7$ & $0.605$ & $p_7$ &  $1$ day$^{-1}$ & &  \\ \hline
\end{tabular}
\label{tab:param1}
\end{center} \end{table}
\\
To make comparison easier, for the continuous maturation model, we replace the interval $[0,x^*]$ by the interval $[1,7]$ (7 being the number of maturing steps in the discrete model) and we  define $a(x)$ and $p(x)$ based on parameters in Table \ref{tab:param1} by piecewise linear continuous functions with values $a_i$ and $p_i$ at $x=i.$ We take them along a regular grid  to obtain approximations of $a(x)$ and $p(x)$.

In Figure \ref{fig:discrete}, the results of the discrete model are identical with the ones of the continuous model if the grid is equal to $X=[1,2,...7]$ (case $I=6$). If the grid becomes finer, we observe a slower convergence toward the steady state together with an increase of the relative importance of the stem cell population. Though unrealistic from a biological viewpoint, it was expected by the derivation of the continuous model from the discrete one. It shows that the analogy between the two models is limited. They exhibit different quantitative properties (see \cite{Anna2} for a study of the discrete model properties), as well as conditions for nontrivial steady state. Moreover, we see that the typical parameter sizes have to be adapted. Indeed, the time evolution is much too slow compared to experimental data. 

%In order to really compare discrete and continuous models, we need here to take $g(x,v)$ defined by $$g(x)=\Delta x 2 (1-\f{a(x)}{1=kv})p(x)$$ with $\Delta x$ the space (maturation) step of the numerical scheme, as noticed previously in Section \ref{subsec:discretetocontinuous}.

Let us now focus on the stability and instability properties, in order to illustrate the theoretical results of Propositions \ref{prop:instab:simplest} and \ref{prop:instab:originalmodel}. 

Figures \ref{fig:instab2:1} and \ref{fig:instab2:2} are an illustration of the instability case stated in Proposition \ref{prop:instab:originalmodel}. Here, we took a maturity interval $[0,X^*]$ with $X^*=50,$ and a proliferation rate $p(x)=p_w + B \chi_{x\geq Y^*}$ with $p_w=30, \, Y^*=20,\,B=50.$ We keep $a(x)$ constant equal to $a_w=0.75$ and $k=1.28x10^{-9}$ as in the discrete case. The maturation speed $g(x)$ is given by Equation \eqref{truedata}.
We see that the destabilization is very slow, and our example is very unrealistic, since the stem cell population level is tiny.

In Figures \ref{fig:instab1:1} and \ref{fig:instab1:2}, we illustrate instability in the case of Proposition \ref{prop:instab:simplest}. We have taken here $X^*=1,$ constant proliferation rate $p(x)=p_w=6$ and maturation rate $g(x)=1$ and $\alpha(v)=(\f{2a_w}{1+kv}-1)p_w$ with $a_w=0.75$ and $k=1.28x10^{-9}.$

\section{Final Remarks  \label{concs}}

In this paper we have developed a structured population model of cell differentiation and self-renewal with a nonlinear regulatory feedback between the level of mature cells and the rate of the maturation process. We showed that perturbations in the regulatory mechanism may lead to the destabilization of the positive steady state, which corresponds to the healthy state of the tissue. In particular, we showed that the regulation of stem cells self-renewal is not sufficient for stability of the system and the lack of the regulation on the level of progenitor cells may lead to the persistent oscillations. This and other stability results suggest how imbalanced regulation of cell self-renewal and differentiation may lead to the destabilization of the system, which is observed during development of some cancers, such as leukemias. 
The model developed in this paper is rather general and, after adjusting it to specific biological assumptions, may serve as a tool to explore the role of different regulatory mechanisms in the normal and pathological development.

Comparing the model to its discrete counterpart we addressed the question of the choice of the right class of models, for example discrete compartments versus continuous maturation, punctuated by division events. We showed that the models may exhibit different dynamics. Interestingly, the structure of steady states varies and the discrete compartmental model admits semi-trivial steady states of the form $(0,..,0,\bar u_i,.., \bar u_n)$, which do not exist in the continuous differentiation model. 

To understand the difference between the two models, we derived a limit equation for the discrete model assuming that a continuum of differentiation stages can be defined. The rationale for such assumption is provided by the fact that differentiation is controlled by intracellular biochemical processes, which are indeed continuous in time, at least when averaged over a large number of cells. Consequently, for the proper time scaling we have to assume that commitment and maturation of cell progenitors do not proceed by the division clock (one division = one step in the maturation process) but is a continuous process and can take place between the divisions. This observation explains the fundamental difference between the two models. The structured population model is indeed not a limit of the discrete model with the transitions between compartments correlated to the division of the cells. However, the models can exhibit exactly the same dynamics for a suitable choice of the maturation rate function $g$. 
%Marie nov 29
%The difference between the qualitative behavior of the solutions of structured population model and the discrete model indicates also that %the numerical discretization of the model may alter its dynamics and analytical understanding of the model dynamics is an important issue.

\begin{figure}[htbp] 
\begin{center}
\begin{minipage}{17cm}
\includegraphics[width=8cm, height=7cm]{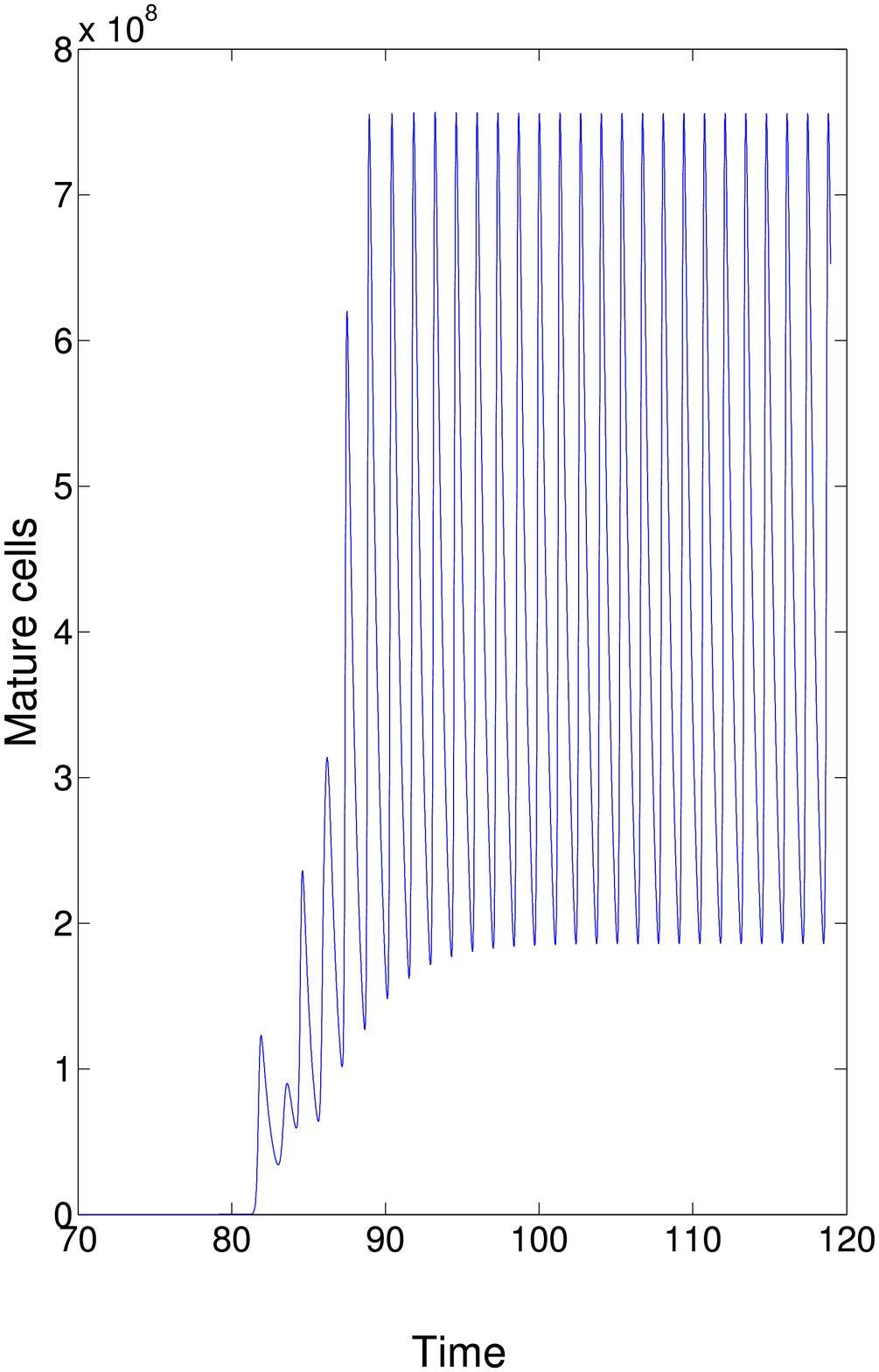} \quad \includegraphics[width=8cm,height=7cm]{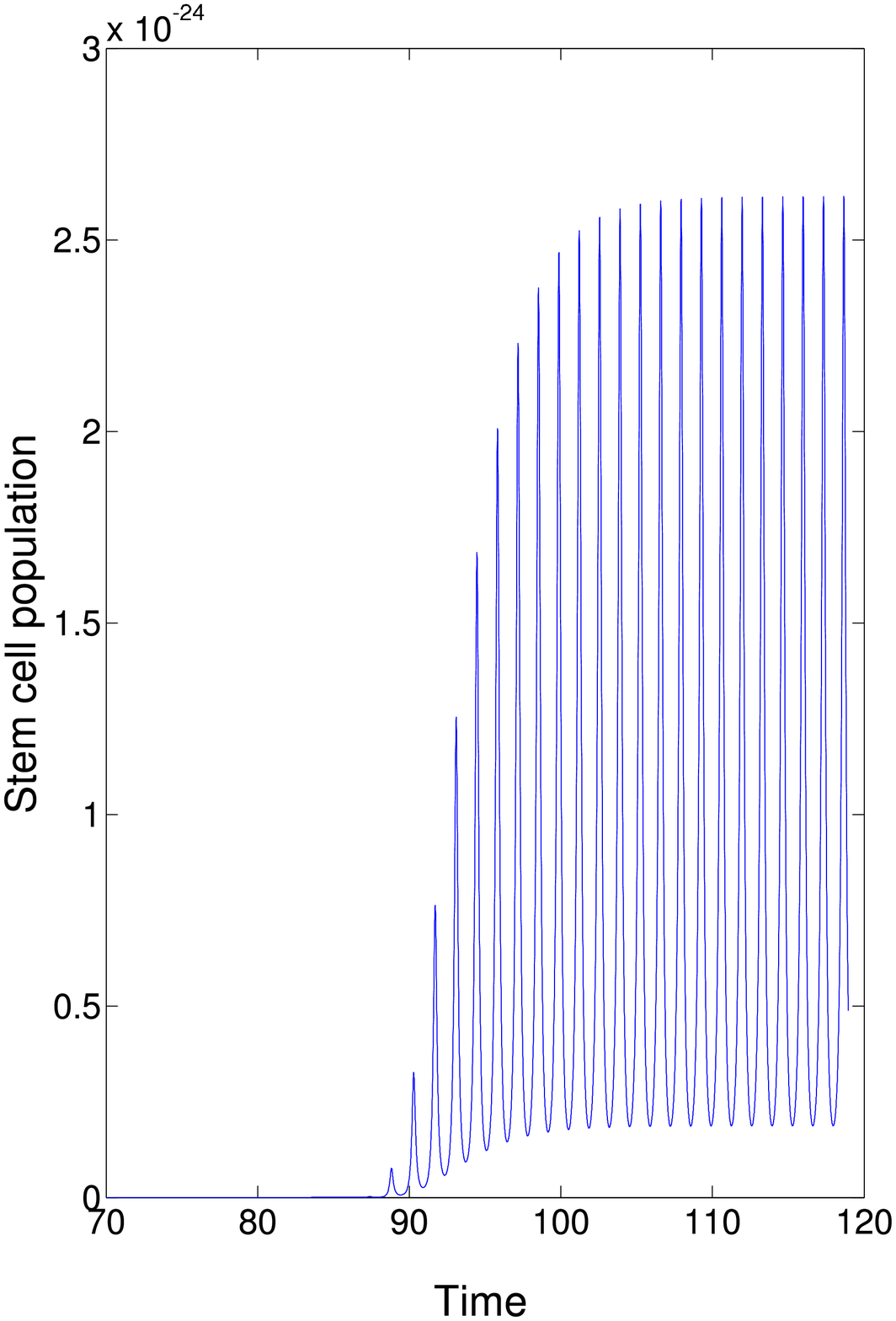}\end{minipage} \end{center}\vspace{-0.5cm}
\caption{\label{fig:instab2:1} Example of instability, in illustration of Proposition \ref{prop:instab:originalmodel}. Left: evolution of mature cells. Right: evolution of stem cells.}
\end{figure}

\begin{figure}[htbp] 
\begin{center}
\begin{minipage}{17cm}
\includegraphics[width=8cm, height=7cm]{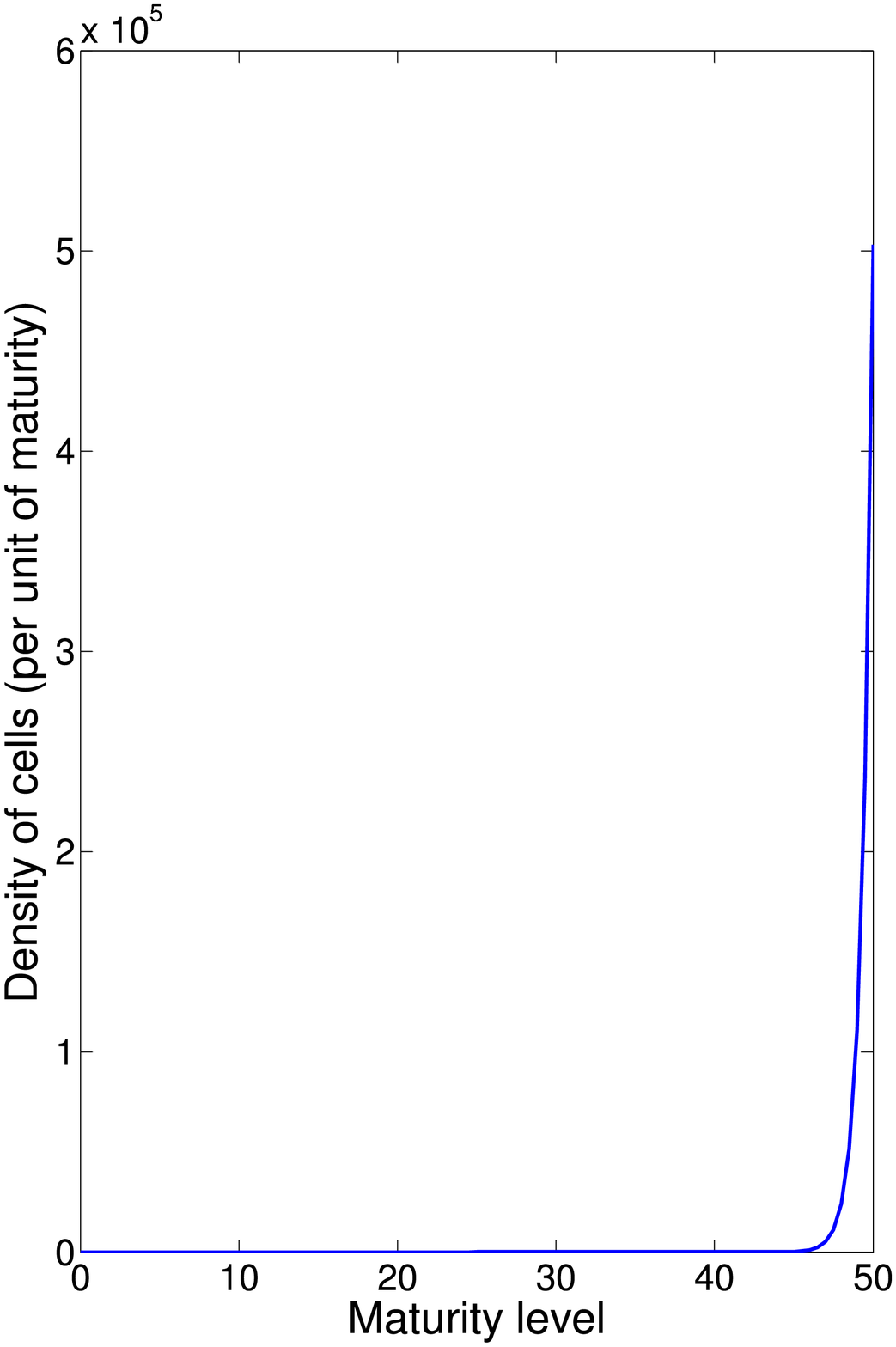} \quad \includegraphics[width=8cm,height=7cm]{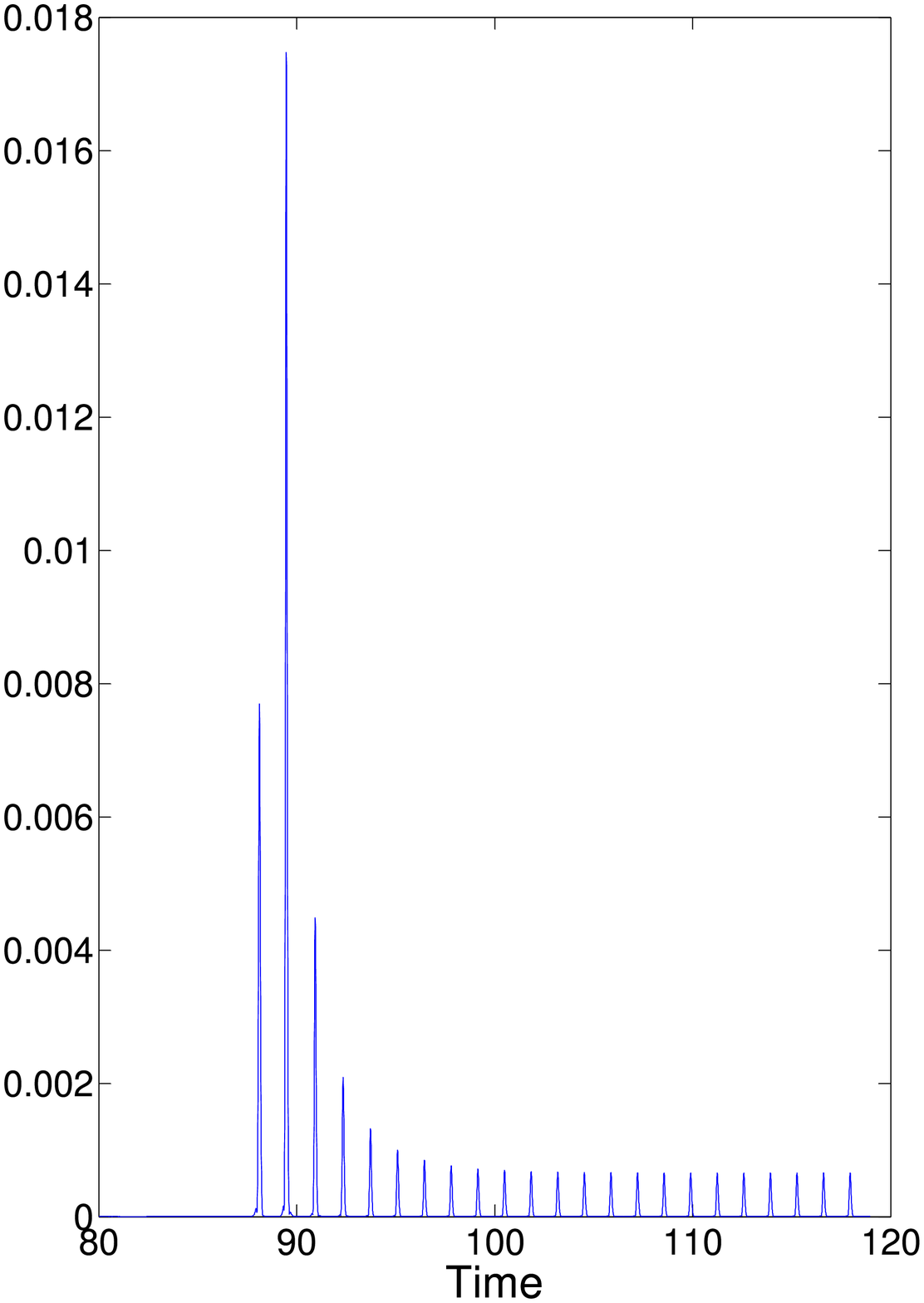}\end{minipage} \end{center}\vspace{-0.5cm}
\caption{\label{fig:instab2:2} Same case as in Figure \ref{fig:instab2:1}. Left: final distribution of cells according to their maturity level. Right: time evolution of $\sqrt{\f{\int |\f{\p}{\p t} u(x,t)|^2 dx}{\int |u(x,t)|^2 dx}}$, to measure the trend to a stable maturity level distribution.}
\end{figure}


\newpage 
\section*{Acknowledgments}
\jpz{}{
The authors were supported by the CNPq-INRIA agreement INVEBIO. 
JPZ was supported by CNPq under grants 302161/2003-1 and
474085/2003-1. AM-C  was supported by ERC Starting Grant "Biostruct", Emmy Noether Programme of German Research Council (DFG) and WIN Kolleg of Heidelberg Academy of Sciences and Humanities. JPZ and BP are thankful to the RICAM Special Semester on Quantitative Biology Analyzed by Mathematical Methods, October 1st, 2007
-January 27th, 2008, organized by RICAM, Austrian Academy of Sciences, to the MathAmSud Programm NAPDE 
and to the International Cooperation Agreement Brazil-France.  
}

\newpage

\appendix

\section{Appendix: Proofs of the Results in  Section~\ref{sec:linearised:original}}

\subsection{The Characteristic Equation in the Case Derived from the Discrete Model}
Using definitions given in \eqref{truedata}, Equation \eqref{eq:eigenvalue:general} can be rewritten as
\begin{equation}\label{eq:eigenvalue:alpha}
\lambda + \mu- \f{d g}{d v} (x^{*},\bar v) \bar u(x^{*})   =  \biggl(-\f{k}{2a_w} p_w  \f{\bar g(x^*)}{\lambda} \bar u(x^*)   +  e^{\int_{0}^{x^*} \frac{p(\xi)}{\bar g(\xi)}d\xi}  \int\limits_0^{x^*} f(s)e^{\int\limits_0^s \f{\lambda }{g(\sigma,\bar v)}d\sigma}ds\biggr) e^{-\int\limits_0^{x^*} \f{\lambda}{g(s,\bar v)}ds}.
\end{equation}
Calculating the term with $f$ we obtain 
$$\int\limits_0^{x^*} f(x)e^{\int\limits_0^x \f{\lambda }{g(s,\bar v)}ds}dx=
\int\limits_0^{x^*} - \partial_x[\f{\partial g}{\partial v} (x,\bar v) \bar u(x)] e^{\int\limits_0^x \f{\lambda -p(s) }{g(s,\bar v)}ds}dx$$
$$=+\int\limits_0^{x^*}  \f{\partial g}{\partial v} (x,\bar v) \bar u(x) \partial_x[e^{\int\limits_0^x \f{\lambda -p(s)}{g(s,\bar v)} ds}]dx - \f{\partial g}{\partial v} (x^*,\bar v) \bar u(x^*)e^{\int\limits_0^{x^*} \f{\lambda -p(s)}{g(\sigma,\bar v)}ds}+\f{\partial g}{\partial v} (0,\bar v) \bar u(0)$$
$$=+\int\limits_0^{x^*}  \f{\partial g}{\partial v} (x,\bar v) \bar u(x) \f{\lambda -p(x)}{g(x,\bar v)}e^{\int\limits_0^x \f{\lambda -p(s)}{g(s,\bar v)} ds}dx - \f{\partial g}{\partial v} (x^*,\bar v) \bar u(x^*)e^{\int\limits_0^{x^*} \f{\lambda -p(s)}{g(\sigma,\bar v)}ds}+\f{k}{2a_w}p_w \bar w.$$
Inserting it into Equation \eqref{eq:eigenvalue:alpha} leads to
$$\lambda + \mu =  \biggl(-\f{k}{2a_w} p_w  \f{\bar g(x^*)}{\lambda} \bar u(x^*)   +  e^{\int_{0}^{x^*} \frac{p(\xi)}{\bar g(\xi)}d\xi} \bigl( 
\int\limits_0^{x^*}  \f{\partial g}{\partial v} (x,\bar v) \bar u(x) \f{\lambda -p(x)}{\bar g(x)}e^{\int\limits_0^x \f{\lambda -p(s)}{\bar g(s)} ds}dx + \f{k}{2a_w}p_w \bar w\bigr)
\biggr) e^{-\int\limits_0^{x^*} \f{\lambda}{\bar g(s)}ds}.
$$
Due to the definition of $\bar w$ given by Equation \eqref{StationaryW}, the first and the last term of the right-hand side can be written together as
$$\lambda + \mu =  \biggl(\f{k}{2a_w}  \bar g(x^*) \bar u(x^*) \bigl(-\f{p_w}{\lambda}+1\bigr)   +  e^{\int_{0}^{x^*} \frac{p(\xi)}{\bar g(\xi)}d\xi} \bigl( 
\int\limits_0^{x^*}  \f{\partial g}{\partial v} (x,\bar v) \bar u(x) \f{\lambda -p(x)}{\bar g(x)}e^{\int\limits_0^x \f{\lambda -p(s)}{\bar g(s)} ds}dx \bigr)
\biggr) e^{-\int\limits_0^{x^*} \f{\lambda}{\bar g(s)}ds}.
$$
Using Equation (\ref{StationaryU}) we obtain
$$\bar u(x) \f{\lambda -p(x)}{\bar g(x)}e^{\int\limits_0^x \f{\lambda -p(s)}{\bar g(s)} ds}
=\bar g (x^*)\bar u (x^*)e^{-\int\limits_0^{x^*} \f{p(s)}{\bar g(s)} ds}\f{\lambda -p(x)}{\bar g(x)^2}e^{\int\limits_0^x \f{\lambda}{\bar g(s)} ds}.
$$
So finally, using Equations \eqref{SteadyState3} and \eqref{StationaryV} and taking $\bar g(x^*) \bar u(x^*)=\f{\mu}{k}(2a_w-1)$, we obtain
the expression given by Equation \eqref{eq:eigenvalue:alpha2}.
\subsection*{Proof of Proposition \ref{prop:instab:originalmodel}}

The local linear stability is equivalent to the fact that all eigenvalues $\lambda \in {\mathbb{C}}$, given by  solutions of Equation \eqref{eq:eigenvalue:alpha3}, have negative real parts.\\
{\bf First step:} $\lambda$ is a solution of the following equation
\begin{equation}\label{eq:eigenvalue:alpha4}\lambda^2 + C\lambda + D=- \f{\mu}{2a_w}(2a_w-1) \lambda \int\limits_0^{x^*}  b(x) e^{\f{\lambda }{p_w} (x-x^*)}dx 
\end{equation}
with $C=\f{\mu}{2a_w} >0,$ $D=p_w \mu \f{2a_w-1}{2a_w} >0,$ and $b(x)=\f{p(x)-p_w}{p_w} \geq 0$ a non-decreasing function. Indeed, with this definition of $b(x)$ we rewrite Equation \eqref{eq:eigenvalue:alpha3} in the form
$$\lambda + \f{\mu}{2a_w} =- \f{\mu}{2a_w}(2a_w-1)\biggl(  \f{p_w}{\lambda}   +    
\int\limits_0^{x^*}  (1+b(x))e^{\f{\lambda }{p_w} x}dx 
\biggr) e^{-\f{\lambda}{p_w}x^*}.$$
Integrating by parts we obtain
$\f{p_w}{\lambda}   +    \int\limits_0^{x^*}  e^{\f{\lambda }{p_w} x}dx=\f{p_w}{\lambda} e^{\f{\lambda }{p_w} x^*}.$ \\
{\bf Second step:} the limiting case is for $b(x)=0,$ \emph{i.e.}, $p$ independent of $x.$ In this case, the eigenvalues are given by 
$$\lambda_\pm=\f{-C \pm \sqrt{C^2 -4D}}{2}.$$
For $C^2-4D >0$ these two eigenvalues are negative. If $C^2-4D<0,$ they are complex conjugated with negative real parts. 
In any case, the steady state is locally linearly stable and the first part of the proposition is proved.\\
{\bf Third step:} In the general case, in order to study the sign of the real part of the eigenvalues $\lambda$, we look for values of the parameters such that $\lambda=i\omega$ with $\omega\in\R.$ It corresponds to a Hopf bifurcation and it leads to:
$$-\omega^2 + iC\omega + D=- \f{\mu}{2a_w}(2a_w-1) i\omega \int\limits_0^{x^*}  b(x) \cos(\f{\omega }{p_w} (x-x^*)) dx + \f{\mu}{2a_w}(2a_w-1) \omega \int\limits_0^{x^*}  b(x)\sin(\f{\omega}{p_w} (x-x^*)) dx.
$$
Taking the imaginary part of this equation yields, since $\omega \neq 0:$
$$C = - \f{\mu}{2a_w}(2a_w-1)  \int\limits_0^{x^*}  b(x) \cos(\f{\omega }{p_w} (x-x^*)) dx
=- \f{\mu}{2a_w}(2a_w-1)  \int\limits_0^{x^*}  b(x^*-y) \cos(\f{\omega }{p_w} y) dy
$$
Since $b$ is increasing, $b(x^*-\cdot)$ is decreasing. This leads to
$$\int\limits_0^{Min(x^*,\f{p_w}{\omega}\pi)}  b(x^*-y) \cos(\f{\omega }{p_w} y) dy \geq 0.$$
Indeed, either $x^* \leq \f{\pi}{2}\f{p_w}{\omega},$ in which case it is evident because $b(x^*-y) \cos(\f{\omega }{p_w} y)  \geq 0$ for all $0\leq x\leq x^*,$ or $x^* \geq  \f{\pi}{2}\f{p_w}{\omega}$ and we can write
$$\int\limits_0^{Min(x^*,\f{p_w}{\omega}\pi)}  b(x^*-y) \cos(\f{\omega }{p_w} y) dy \geq b(x^*-\f{\pi}{2}\f{p_w}{\omega}) \int\limits_0^{\f{\pi}{2}\f{p_w}{\omega}}  \cos(\f{\omega }{p_w} y) dy - b(x^*-\f{\pi}{2}\f{p_w}{\omega}) |\int\limits_{\f{\pi}{2}\f{p_w}{\omega}}^{\f{p_w}{\omega}\pi} \cos(\f{\omega }{p_w} y) dy  = 0.$$
To end the proof, let us simply exhibit an example where instability can occur: Let $\chi$ denote the Heaviside function. Defining
$$b(x)=B\chi_{y^*\leq x \leq x^*},$$
we compute explicitely
$$-\omega^2 + i\f{\mu\omega}{2a_w} + p_w\mu \f{2a_w-1}{2a_w}=- p_w{\mu}\f{2a_w-1}{2a_w} i B  \sin(\f{\omega }{p_w} (x^*-y^*)) + p_w\mu\f{2a_w-1}{2a_w} B    \bigl( \cos(\f{\omega }{p_w} (x^*-y^*))-1\bigr).
$$
It provides two relations
$$\begin{array}{l}
\f{\omega}{p_w}=- (2a_w-1) B  \sin(\f{\omega }{p_w} (x^*-y^*)),
\quad 
\f{\omega^2}{p_w^2}\f{2a_w}{2a_w-1}  = \f{\mu}{p_w} \biggl(1+ B    \bigl(1- \cos(\f{\omega }{p_w} (x^*-y^*))\bigr)\biggr).\\ \\
\end{array}
$$
We can see that there exist many sets of parameters such that both relations are satisfied. Given $a_w,$ $B,$ $p_w,$ and $\omega$  such that 
$$\f{\omega}{p_w(2a_w-1)B}\leq 1,
$$ 
we can always find $x^*-y^*$ such that the first relation is satisfied, and then fix $\mu$ using the second relation.
\qed

%-------------------------------------------------------

\begin{thebibliography}{99}

\bibitem{Adimy1} M.~Adimy, F.~Crauste and L.~Pujo-Menjouet, \textit{On the stability of a maturity structured model of cellular proliferation,} Dis. Cont. Dyn. Sys. Ser. A, 12 (3): 501--522, 2005.
\bibitem{Adimy2} M.~Adimy, F.~Crauste and S.~Ruan, \textit{A mathematical study of the hematopoiesis process with applications to chronic myelogenous leukemia}, SIAM J. Appl. Math., 65 (4): 1328--1352, 2005. 
\bibitem{Al-Hajj} M.~Al-Hajj and M.F.~Clarke, \textit{Self-renewal and solid tumor stem cells,} Oncogene 23: 7274-7282, 2004.
\bibitem{Beachy} P.A.~Beachy, S.S.~Karhadkar and D.M.~Berman,  \textit{Tissue repair and stem cell renewal in carcinogenesis,} Nature 432: 324-331, 2004.
 \bibitem{2Mackey} J.~Belair, M.C.~Mackey and J.~Mahaffy, \textit{Age structured and two delay models for erythropoiesis}, Math. Biosci. 128:  317--346, 1995. 
\bibitem{Mackey} C.~Colijn and M.C.~Mackey, \textit{A mathematical model of hematopoiesis--I. Periodic chronic myelogenous leukemia}, J.~Theor.~Biol. 237: 117--132, 2005. 
\bibitem{GV} J.-F.~Collet, T.~Goudon, F.~Poupaud and A.~Vasseur,
{\it The {B}ecker--{D}\"{o}ring System and Its {L}ifshitz--{S}lyozov Limit}, {SIAM J. on Appl. Math.} 62 (5): 1488--1500, 2002.
\bibitem{Diekmann} O.~Diekmann, S.A~van Gils, S.M.~Verduyn-Lunel and H.C.~Walther, \emph{Delay equations, Functional-, Complex and Nonlinear Analysis}, ApMS, 110, NY, 1995.
\bibitem{Dontu} G.~Dontu, M.~Al-Hajj, W.M.~Abdallah, M.F.~Clarke and M.S.~Wicha, 
  \emph{Stem cells in normal breast development and breast cancer},
Cell Prolif. 36: 59-72, 2003. 
\bibitem{DGL} M.~Doumic, L.~Goudon and T.~Lepoutre, \emph{{S}caling limit of a discrete prion dynamics model},
{Comm. in Math. Sc.} 7(4): 839--865, {2009}.
 \bibitem{Kim} M.~Doumic, P.~Kim, B.~Perthame   \emph{Stability Analysis of a Simplified Yet Complete Model for Chronic Myelogenous Leukemia}, Bull. of Math. Biol. (doi: 10.1007/s11538-009-9500-0).
 \bibitem{He} S.~He, D.~Nakada, S.J.~Morrison,  \emph{Mechanisms of stem cell self-renewal}, Annu. Rev. Cell. Dev. Biol. 25: 377-406, 2009.
\bibitem{Hope} K.J.~Hope, L.~Jin and J.E.~Dick, \emph{Acute myeloid leukemia originates from a hierarchy of leukemic stem cell classes that differ in self-renewal capacity}. Nat Immunol. 5: 738-743, 2004. 
\bibitem{Klaus} J.~Klaus, D.~Herrmann, I.~Breitkreutz, U.~Hegenbart, U.~Mazitschek, G.~Egerer,
F.W.~Cremer, R.M.~Lowenthal, J.~Huesing, S.~Fruehauf, T.~Moehler, A.D.~Ho,and
H.~Goldschmidt, \emph{Effect of cd34 cell dose on hematopoietic reconstitution and
outcome in 508 patients with multiple myeloma undergoing autologous peripheral
blood stem cell transplantation,} Eur J Haematol. 78(1): 21--28, 2007.
\bibitem{Lord} B.I.~Lord, \emph{Biology of the haemopoietic stem cell}, In C. S. Potten, editor, Stem cells: 401--422. Academic Press, Cambridge, 1997.
\bibitem{MagalRuan} P.~Magal and S.~Ruan, \emph{Center Manifolds for Semilinear
Equations with Non-dense Domain and Applications on Hopf Bifurcation in Age
Structured Models}, Memoirs of the AMS. 202(951), 2009.
\bibitem{Anna1} A.~Marciniak-Czochra, T.~Stiehl, W.~J\"ager, A.~Ho and W.~Wagner, \emph{ {M}odeling of
asymmetric cell division in hematopoietic stem cells - regulation of self-renewal
is essential for efficient repopulation,} Stem Cells Dev. 18(3): 377-385, 2009.
\bibitem{Anna2} A.~Marciniak-Czochra, T.~Stiehl,  \emph{{C}haracterization of stem cells using mathematical
models of multistage cell lineages,} {Math. Comp. Models.} (doi: 10.1016/ j.mcm.2010.03.057).
\bibitem{Moore} K.A.~Moore and I.R.~Lemischka \emph{Stem cells and their niches,} Science 311: 1880-1885, 2006.
\bibitem{Morrison}
S.J.~Morrison and J.~Kimble, \emph{ Asymmetric and symmetric stem cell divisions in development and cancer}, Nature 441: 1068-1074, 2006. 
\bibitem{Weissmann1} T.~Reya, S.J.~Morrison, M.F.~Clarke and I.L.~Weissman, \emph{Stem cells, cancer, cancer stem cells}, Nature 414: 105-11 2001.
\bibitem{Jiang} 
G.~Shan Jiang and D.~Peng, \emph{Weighted eno schemes for Hamilton-Jacobi equations}, SIAM J. Sci. Comput, 21: 2126-2143, 1997.
\bibitem{Schroeder} T.~Schroeder \emph{Asymmetric Cell Division in Normal and Malignant Hematopoietic Precursor Cells} Cell Stem Cell 1: 479-481, 2007.
\bibitem{Shu} C.-W.~Shu, \emph{Essentially non-oscillatory and weighted essentially non-oscillatory schemes for hyperbolic
conservation laws,} Springer: 325-432, 1998.
\bibitem{ThomasDiploma}T.~Stiehl. Diploma thesis, University of Heidelberg, 2009.
\bibitem{Thornley} I.~Thornley, D.R.~Sutherland, R.~Nayar, L.~Sung, M.H.~Freedman, and H.A.~Messner,
\emph{Replicative stress after allogeneic bone marrow transplantation: changes
in cycling of cd34+cd90+ and cd34+cd90- hematopoietic progenitors,} Blood,
97(6): 1876–8, 2001.
\bibitem{Uchida} N.~Uchida, W.H.~Fleming, E.J.~Alpern, and I.L.~Weissman, \emph{Heterogeneity of
hematopoietic stem cells,} Curr. Opin. Immunol., 5(2): 177–184, 1993.
\bibitem{Weissman} I.L.~Weissman, \emph{Stem cells: units of development, units of regeneration, and units in evolution},
Cell 100(1): 157-68, 2000. PMID: 10647940.
\end{thebibliography}
\end{document}